\documentclass[twoside,reqno,draft]{amsart}
\usepackage{amssymb}
\usepackage[dvips]{graphics}
\newtheorem{thm}{Theorem}[section]
\newtheorem{lem}[thm]{Lemma}

\newtheorem{cor}[thm]{Corollary}

\newtheorem{prop}[thm]{Proposition}

\renewcommand{\qed}{\hfill$\Box$ \par \medskip}

\newcommand{\ignore}[1]{}

\newcommand{\Vol}{\mathrm{Vol}\,}

\newcommand{\ball}[2]{B_{#1}(#2)}

\newcommand{\dist}{\mathrm{dist}}

\newcommand{\eqdef}{:=}

\newcommand{\norm}[2][{}]{\left\|#2\right\|_{#1}}

\newcommand{\origin}{\mathbf{0}}

\newcommand{\pair}[2]{\left\langle #1,#2 \right\rangle}

\newcommand{\R}{\mathbf{R}}

\newcommand{\N}{\mathbf{N}}

\newcommand{\supp}{\mathrm{supp}\,}

\newcommand{\cC}{\mathcal{C}}

\newcommand{\cE}{\mathcal{E}}

\title{Energy and Invariant Measures for Birational Surface Maps}

\author{Eric Bedford \& Jeffrey Diller}
\address{Department of Mathematics\\
Indiana University\\
Bloomington, IN 47405}
\email{}
\address{Department of Mathematics\\
         University of Notre Dame\\
         Notre Dame, IN 46656}
\email{diller.1@nd.edu}


\begin{document}
\maketitle

\section*{0.  Introduction}
\noindent Let $X$ be a compact K\"ahler surface, and let
$f:X\to X$ be a bimeromorphic mapping.  We consider $(f,X)$ as
a dynamical system, which means that we consider the behavior of
the iterates $f^n=f\circ\cdots\circ f$ as $n$ tends to
infinity.  Since $f$ is invertible we may consider both forward
and backward dynamics, i.e., $f^n$ as $n\to+\infty$ and
$n\to-\infty$.  A meromorphic map of a surface is holomorphic
outside a finite set, $I(f)$,  of points which are blown up
to curves.  Thus $f$ is not in general a continuous
map, so it is not clear to what extent there is a standard
category of dynamical systems into which such an object falls.

We consider two bimeromorphic maps to be equivalent
if they are bimeromorphically conjugate.  Two complex surfaces
can be bimeromorphically equivalent, however, without being
homeomorphic.  One approach that has proved fruitful in  complex
dynamics is to start with the induced action $f^*$ on the
cohomology group $H^{1,1}(X)$.  A question that arises when
$f$ has points of indeterminacy is whether the passage to
cohomology is natural for the dynamics, i.e., whether
$(f^n)^*=(f^*)^n$.  This happens exactly
when the condition
\begin{equation}\bigcup_{n\ge0}f^{-n}I(f)\cap\bigcup_{n\ge0}f^nI(f^{-1})
=\emptyset
\label{eq1}\end{equation}
holds.  This condition may be viewed as a separation between the
obstructions to forward and backward dynamics.  Diller and Favre [DF]
showed that any bimeromorphic surface map
$f:X\to X$ is bimeromorphically equivalent to a map $\hat f:\hat
X\to\hat X$ for which (\ref{eq1}) holds.  In general, the spectral radius
$\rho$ of $\hat f^*$ on
$H^{1,1}(X)$ is greater than or equal to 1, and it was shown in [DF] that
if
$\rho=1$, then either $f$  is a dynamically trivial automorphism,
or $f$ preserves a rational or elliptic fibration and exhibits a
dynamic which is essentially  one-dimensional.

We assume in this paper that
$\rho>1$.  In this case there are stable/unstable
currents
$\mu^\pm$ whose cohomology classes generate the $f^*$ and $f_*$
eigenspaces for $\rho$, and in fact
$f^*\mu^+=\rho\mu^+$ and
$f_*\mu^-=\rho\mu^-$.
The currents $\mu^\pm$ carry geometric information of
(complex) dimension 1 and are useful in analyzing the dynamics
of
$f$.

A natural hope is that the wedge product
$\mu:=\mu^+\wedge\mu^-$ might define an invariant measure that serves
as a bridge between the action of $f^*$ on $H^{1,1}$ and the ergodic
properties of $f$ on $X$.  This was shown to happen
for polynomial automorphisms of
${\bf C}^2$ in the papers [BS] and [FS]; for automorphisms of K3
surfaces in [C]; and for
certain birational maps in [D] and [G].  Typically one considers the
positive, closed currents
$\mu^\pm= dd^cg^\pm$ in terms of local potentials.  The operation of
wedge product is then interpreted in terms of the so-called
complex Monge-Amp\`ere operator $dd^cg^+\wedge dd^cg^-$.
As is well known, this operation is possible if at least one of
the potentials $g^+$ or $g^-$ is locally bounded.  And this is
what happens in all of the papers cited above.  On the other
hand, it is possible for both local potentials $g^+$ and $g^-$ to
be locally unbounded at a point, as is the case for the ``golden
mean'' family, which was analyzed in detail in [BD].

The condition
\begin{equation}
\overline{\bigcup_{n\ge0}f^{-n}I(f)}
\cap\overline{\bigcup_{n\ge0}f^nI(f^{-1})}
=\emptyset
\label{eq2}\end{equation}
was introduced in [Dil2], and it was shown to be equivalent to
the condition that for each point there is a neighborhood on which one
of the local potentials $g^+$ or $g^-$ is continuous.   In this paper we
employ a quantitative condition stronger than (\ref{eq1}) and weaker than
(\ref{eq2}):
\begin{equation}
\sum_{n\ge0} \rho^{-n}\log dist(f^nI(f^{-1}),I(f))>-\infty.
\label{eq3}\end{equation}
By Theorem \ref{mu15} this is equivalent to $g^+(x)>-\infty$ for all
$x\in I(f)$.

\newtheorem*{thm1}{Theorem}
\begin{thm1}  If (\ref{eq3}) holds, then $\mu:=\mu^+\wedge\mu^-$ is
a probability measure that puts no mass on any algebraic set;
$\mu$ is invariant and mixing for $f$.  Further,
\begin{equation}
\int\left|\log||Df||^{}\,\right|\, \mu<\infty,
\label{eq3a}
\end{equation}
and thus the Lyapunov exponents of $f$ with respect to $\mu$ are
well-defined and finite.  Finally,  the
Lyapunov exponents satisfy
$$\lambda^-\le-{\log\rho\over8}<0<
{\log\rho\over 8}\le \lambda^+,$$ and thus
$\mu$ is a hyperbolic measure of saddle type.
\end{thm1}
The finiteness of the integral in (\ref{eq3a}), and thus the finiteness
of the Lyapunov exponents, seems to be closely linked with condition
(\ref{eq3}).  On the other hand, Favre [Fav3]  has constructed a mapping
which does not satisfy (\ref{eq3}).  Favre's example depends on the
existence of an invariant complex line whose rotation number satisfies a
delicate number-theoretic property.

We define $\mu^+\wedge\mu^-$ using an ``energy''
approach to interpret the complex Monge-Amp\`ere operator, as was done in
[BT] and [B].  Specifically, if $T$ is a positive, closed
(1,1)-current, then we define the energy of a function $\varphi$ to
be
$$\cE_T(\varphi):=\int d\varphi\wedge d^c\varphi\wedge T.$$
The approach from [BT] is that if $\varphi$ is
essentially psh, and if $\cE_T(\varphi)<\infty$, then
$dd^c\varphi\wedge T$ defines a measure, and $\varphi$ is
integrable with respect to this measure.  In the situation at
hand, we will show that
$\cE_T(\varphi)<\infty$ for
$T=\mu^-$ and
$\varphi=g^+$.

The currents $\mu^\pm$ are obtained dynamically by starting with
K\"ahler forms
$\beta_1$ and
$\beta_2$ and taking  normalized limits of
pullbacks:
\begin{equation}\mu^+=c_1\cdot\lim_{n\to\infty}\rho^{-n}f^{*n}\beta_1,\ \
\
\mu^-=c_2\cdot\lim_{n\to\infty}\rho^{-n}f^n_*\beta_2.  \label{eq3b}
\end{equation}
We show in Corollary \ref{mu3}  that the measure $\mu$ is also
obtained as
\begin{equation}\mu=c\cdot\lim_{n,m\to\infty}\rho^{-n-m}f^{*n}\beta_1\wedge
f^m_*\beta_2. \label{eq3c}
\end{equation}

The contents of the paper are as follows.  In \S1 we discuss the
pullbacks of currents and the associated (local) potential functions.
The fact of convergence in (\ref{eq3b}) was established in [DF].
However, in order to pass from (\ref{eq3b}) to (\ref{eq3c}), we need to
know how the intermediate pullbacks depend on
$\beta_1$ and $\beta_2$.  This dependence is clarified in \S2.  In \S3
we discuss properties of the energy integral.  In \S4 we discuss
condition (\ref{eq3}); we show that when (\ref{eq3}) holds the gradients
of the local potentials of $\mu^\pm$ belong to $L^2$.  Thus
$\mu:=\mu^+\wedge\mu^-$ is well defined.  We show in Theorem
\ref{invariance} that
$\mu$ is invariant.  \S5 is dedicated to showing that $\mu$ is mixing,
and
\S6 gives the estimates on the Lyapunov exponents.

\section{Pullbacks under Birational Maps}
\noindent Throughout this paper we let $X$ denote a compact K\"ahler surface
endowed with the  hermitian
metric associated to a fixed K\"ahler form $\beta$.  Let $f:X\to X$ be
a  bimeromorphic
self-map.  That is, there
is a compact surface $\Gamma$ (the desingularized graph of $f$) with
proper modifications (i.e. generically injective holomorphic maps)
$\pi_1,\pi_2:\Gamma\to X$ such that $f=\pi_2\circ\pi_1^{-1}$.  The set
$$\cC(\pi_j)\eqdef\{x\in\Gamma:\#(\pi_j^{-1}\pi_j(x))>1\} = \{x\in\Gamma:
\text{dim}(\pi_j^{-1}\pi_j(x))=1\}$$ is the critical  set for
$\pi_j$.
The images
$I(f)\eqdef\pi_1(\cC(\pi_1))$ and $\cC(f) \eqdef\pi_1(\cC(\pi_2))$ are
the \emph{indeterminacy} and \emph{critical} sets, respectively, for $f$.
Note that in this case the critical set is actually an exceptional set,
since the irreducible components are mapped to points. It is shown in
[DF], section 2 that after a finite number of blow-ups we may assume that
(1) holds.
In this case, $I(f^n)=\bigcup_{j=0}^{n-1}f^{-j}I(f).$

Since $f$ is ill-defined at points of indeterminacy, it is useful to
adopt some conventions concerning images of points and curves under
$f$.  Given any $x\in X$, we set $f(x) = \pi_2(\pi_1^{-1}(x))$ with
the effect that $f(x)$ is a point if $x\notin I(f)$ and a component of
$\cC(f^{-1})$ otherwise.  Given
any curve $V\subset X$, we set $f(V) = \overline{f(V\setminus I(f))}$.
For irreducible $V$, it follows that $f(V)$ is a point if $V\subset\cC(f)$
and an irreducible curve if not.

\begin{prop}
\label{derivative}
There exist constants $A,B>0$ such that
$$
\norm{D_xf} \leq A\,\dist(x,I)^{-B}
$$
for all $x\in X$.  Further, given a second point $y\in X$, one has
$$
\dist(f(x),f(y)) \leq A\,\dist(\{x,y\},I)^{-B}\dist(x,y).
$$
\end{prop}

\proof
Choose any hermitian metric on the graph $\Gamma$ of $f$.  Then
$\norm{D\pi_2}$ is uniformly bounded on $\Gamma$, so it suffices to
prove the first inequality for $\pi_1^{-1}$ in place of $f$.

In local coordinates, the entries of $D\pi_1$ are holomorphic functions,
so the entries of $(D\pi_1)^{-1}$ are meromorphic functions with poles
in $\cC(\pi_1)$.  Since $\Gamma$ is compact, there are constants
$A, B > 0$ such that
$$
\norm{(D_y\pi_{1})^{-1}} \leq A\,\dist(y,\cC(\pi_1))^{-B}
$$
for all $y\in \Gamma$.  But $\norm{D\pi_1}$
is uniformly bounded on $\Gamma$, so this implies
$$
\norm{D_{\pi_1(y)}(\pi^{-1}_1)} = \norm{(D_y\pi_{1})^{-1}}
\leq A\,\dist(\pi_1(y),I)^{-B}.
$$
The first inequality now follows because $\pi_1$ is surjective.

The second inequality follows from the first by integrating along a path
from $x$ to $y$.
\qed

We consider the hermitian inner product on the set of smooth (1,1)-forms
given by
$$\langle\alpha,\beta\rangle:= \int_X\alpha\wedge\bar\beta.$$
It follows that any smooth (1,1)-form defines an element of the dual space
of (1,1)-forms, and thus defines a (1,1)-current.  The (1,1) cohomology group
$H^{1,1}(X)$ may be given as the smooth, closed (1,1)-forms modulo the exact
ones.  It follows from Stokes' Theorem that the hermitian pairing on
(1,1)-forms induces a pairing on
$H^{1,1}(X)$.  In fact, this pairing is a nondegenerate duality.

If $T$ is a closed (1,1)-current, then $T(d\xi)=0$, which means that
$T$ annihilates all $d$-exact (1,1)-forms.  Thus the restriction of $T$ to
the closed forms defines an element of $H^{1,1}(X)^*$, and there is a
cohomology class $\{T\}\in H^{1,1}(X)$ which represents this restriction in
the sense that $T=\langle\cdot,\{T\}\rangle$.

The ``$\partial\bar\partial$-Lemma'' from K\"ahler geometry (see
[GH, page 149]) asserts that if $T_1$ and $T_2$ are closed (1,1) currents
which define the same cohomology class, then there there is a current $S$ of
degree 0 such that
$$T_1=T_2+dd^cS.$$
  In particular, if $T$ is a closed (1,1)-current on $X$, there is
a smooth (1,1)-form
$\alpha$ defining the cohomology class $\{T\}$, and by the
$\partial\bar\partial$-Lemma, there is a current
$h$ such that
$T=\alpha + dd^ch$.

Next we define the pullback of a smooth form.  If $\alpha$ is a smooth (1,1)
form on
$X$, then
$\pi_2^*\alpha$ is a smooth (1,1) form on $\Gamma$.  By duality,
$\pi^*_2\alpha$ defines a current on
$\Gamma$ of bidegree (1,1).  Thus
$$f^*\alpha:=\pi_{1*}(\pi^*_2\alpha)$$ is a current on $X$.  The pullback
$f^*$ commutes with $d$ and with the complex structure, so closed
(respectively, exact) forms are pulled back to closed (resp.\ exact)
currents of the same bidegree.  This gives a well defined map $f^*$ on
$H^{1,1}(X)$.
  Similarly, we set $f_*\eta\eqdef (f^{-1})^*\eta=\pi_{2*}\pi_1^*\eta$.  In
other words, we set $f_*=(f^{-1})^*$.  Note that
$f^*$ and $f_*$ are adjoint with respect to the intersection form
$\pair{\cdot}{\cdot}$ on cohomology classes, which is to say
$$
\pair{f^*\alpha}{\beta} =
\langle\pi^*_2\alpha,\pi_1^*\beta\rangle =\pair{\alpha}{f_*\beta}.
$$

We can also define the pullback $f^*T$ if $T$ is a positive, closed
(1,1)-current on
$X$.  By pulling back local potentials of $T$, we may define $f^*T$ on
$X-I(f)$.  Now for any
$x\in I(f)$, we may choose a pseudoconvex neighborhood $U$ of $x$ with
$H^2(U-\{x\})=0$.  Thus there is a potential $p$ on $U-\{x\}$ such that
$f^*T=dd^cp$ on
$U-\{x\}$.  Since
$p$ is psh on $U-\{x\}$, it follows that $p$ has a psh extension $\tilde p$
to $U$.  We define $f^*T:=dd^cp$ on $U$.

Let $P(X)$ denote the set of upper semicontinuous
functions
$u$ on $X$ such that $dd^c u \geq -c\beta$ for some $c\in\R$.  Such
functions are locally the sum of a psh function and a smooth function.
(Since $X$ is compact, there are no global psh functions.)  Given a finite
set
$S\subset X$, let
$\tilde P(X,S)$ denote those functions
$u\in  P(X)\cap C^\infty(X\setminus S)$ such that
$$
u(x) \geq A\log\dist(x,S) - B
$$
for some $A,B>0$ and all $x\in X$.

\begin{prop}
\label{prelim24}
Suppose that $S\subset X$ is finite and disjoint from $I(f^{-1})$.  Then
$u\in\tilde P(X,S)$ implies that $u\circ f$ is a difference of functions
in $\tilde P(X,f^{-1}(S)\cup I(f))$.
\end{prop}

\proof
Because $u\in \tilde P(X,S)$, we get
$$
0 \leq dd^c u\circ f + cf^*\beta =
        dd^c ( u\circ f + v ) + \beta' \leq dd^c(u\circ f + v) + c\beta,
$$
where $\beta'$ is a smooth $(1,1)$ form cohomologous to $f^*\beta$,
$v\in\tilde P(X,I(f))$, and $c>0$ is chosen large enough that
$c\beta \geq \beta'$.
Moreover, since $u\in\tilde P(X,S)$, we see from Proposition
\ref{derivative} that for $fx$ near $S$, and therefore uniformly far from
$I(f^{-1})$,
\begin{eqnarray*}
u\circ f(x) & \geq & A\log\dist(f(x),S) - B  \\
    & \geq & A\log\dist(x,f^{-1}(S)) - B + C\log\dist(f(x),I(f^{-1})) \\
    & \geq & A\log\dist(x,f^{-1}(S)) - B.
\end{eqnarray*}
Combining the two displayed inequalities, we see that
$$
u\circ f = (u\circ f + v) - v,
$$
where $v\in \tilde P(X,I(f))$, and
$u\circ f + v \in \tilde P(X,f^{-1}(S)\cup I(f))$.
\qed

The Lelong number of a positive closed $(1,1)$ current $T$ at a
point $x\in X$ is the non-negative number
$$
\nu(T,x) \eqdef
\lim_{r\to 0} \frac{C}{r^2} \int_{\ball{x}{r}} \beta\wedge T.
$$
If $u$ is a local potential for $T$ in a neighborhood of $x$, i.e., if $T =
dd^c u$, then
$$
\nu(T,x)  \eqdef \sup\{t\geq 0: u(y) < t\log\dist(x,y) + O(1)\}
$$
  (see [Dem2, Equation 5.5e]).
We use Proposition \ref{derivative} to gain control over the singularities
of pullbacks of smooth (1,1) forms:
\begin{prop}
\label{prelim1}
Let $\omega$ be a K\"ahler form on $X$ and $\omega'$ be a smooth
form cohomologous to $f^*\omega$.  Then we can write
\begin{equation}
\label{prelim2}
f^*\omega = \omega' + dd^c u,
\end{equation}
where $u$ is smooth and negative on $X\setminus I(f)$ and satisfies
$$
A\log\dist(x,I) - B \leq u(x) \leq A'\log\dist(x,I) + B'
$$
for some constants $A,B,A',B' > 0$ and every $x\in X$.
\end{prop}

\proof
If $\omega'$ and $f^*\omega$ represent the same element of $H^{1,1}(X)$,
then there exists a $u$ satisfying (\ref{prelim2}).
The current $f^*\omega$ is positive, so $u$ is given locally as the
sum of a smooth function $u_1$ and a plurisubharmonic function $u_2$.  In
particular, we can assume that $u$ is negative.  The remaining assertion
in the proposition only concerns some (any) choice of $u_2$ in the
neighborhood of a point $y\in I(f)$.

For each component $V'$ of $\pi_1^{-1}(y)$ and its image $V = \pi_2(V')$,
we have
$$
\int_{V'} \pi_2^*\omega = \int_{V} \omega > 0.
$$
The intersection form on $\pi_1^{-1}(y)$ is negative definite, so we can
choose a non-trivial effective divisor $V'$ supported on
$\pi_1^{-1}(y)$ such that $\pi_2^*\omega + [V']$ is cohomologically
trivial near $V$.  In particular, we can write $\pi_2^*\omega + [V'] = dd^c v$
for some function $v$ defined in a neighborhood $U'$ of $\pi_1^{-1}(y)$ and
smooth off $\pi_1^{-1}(y)$.  Therefore $v\circ\pi_1^{-1}$
is a local potential for
$$
\pi_{1*} (\pi_2^*\omega + [V']) = \pi_{1*}\pi_2^*\omega
$$
on the neighborhood $U = \pi_1(U')$ of $y$.  The singularities of $v$ come
entirely from local potentials for $[V']$.  Hence we can arrange
$$
v(x') \geq A\log\dist(p',\pi_1^{-1}(y))
$$
for some $A>0$ and all $x'\in U$.  Finally, since $\pi_1$ is uniformly
Lipschitz, we obtain after adjusting $A$ that
$$
u_2\eqdef v\circ\pi_1^{-1}(x) \geq A\log\dist(x,y),
$$
which finishes the proof of the lower bound for $u_2$.

To obtain the upper bound for $u_2$, we rely on the push-pull
formula [DF, Theorem 3.3] applied to $\pi_1$.  This gives
$$
\pi_1^* f^*\omega = \pi_1^*\pi_{1*}(\pi_2^*\omega) = \pi_2^*\omega +
       [V']
$$
where $V'$ is an effective divisor with support equal to $\cC(\pi_1)$.
In particular, the Lelong number of of the positive current
$\pi_1^* f^*\omega$ is positive everywhere on $\cC(\pi_1)$.
It follows from [Fav1, Theorem 2] that $f^*\omega$ has a positive
Lelong number at each point in $I(f) = \pi_1(\cC(\pi_1))$.  We conclude
that any local potential $u_2$ for $f^*\omega$ near $y\in I(f)$ must
satisfy
$$
u_2(x) \leq A'\log\dist(x,y) + B'
$$
for some $A',B'>0$.
\qed

\section{Invariant Cohomology Classes and Currents}
\label{prelim}
\noindent The condition
(\ref{eq1}) implies that
$(f^n)^* = (f^*)^n$ on $H^{1,1}(X)$ for every $n\in{\bf Z}$, (see [FS]
and [DF, Theorem 1.14]).  In this case the bimeromorphically invariant
quantity
$$
\rho\eqdef \lim_{n\to\infty}\norm{f^{n*}|_{H^{1,1}}}^{1/n} \geq 1
$$
is the modulus of the  largest eigenvalue of $f^*$ on
$H^{1,1}(X)$.  In this
paper, we assume that
\begin{equation}\rho > 1. \label{eq5}
\end{equation}
An element of $H^{1,1}$ is a K\"ahler class if it
contains a K\"ahler form.  We say that a cohomology class is nef if it is
in the closure of the K\"ahler classes.

\begin{thm}
\label{prelim7} If (\ref{eq5}) holds, then
$\rho$ is the unique (counting multiplicity) eigenvalue of $f^*$ of
modulus larger than one; and the associated eigenspace is generated by a
nef class $\theta^+$.
\end{thm}

Since $f^*$ and $f_*$ are adjoint with respect to
the intersection product, $\rho(f^{-1})=\rho(f)$, and so Theorem
\ref{prelim7} yields a class
$\theta^- = \rho^{-1} f_*\theta^-$.  We scale $\theta^\pm$ and
$\beta$ to achieve
\begin{equation}
\pair{\theta^+}{\theta^-} =
\pair{\theta^+}{\beta} =
\pair{\theta^-}{\beta} = 1,\label{eq6}
\end{equation}
and in this case $\theta^+$ and $\theta^-$ are unique.

We fix K\"ahler forms $\omega_1,\dots,\omega_N$ whose cohomology
classes form a basis a for $H^{1,1}(X)$, and we let $\Omega$ denote
the linear span of these forms.  We also assume for convenience that
$\Omega$ contains the K\"ahler form
$\beta$ corresponding to the metric on $X$.  We endow $\Omega$ with
the norm
$\Vert\omega\Vert=(\sum|c_j|^2)^{1/2}$ where $\omega=\sum c_j\omega_j$.
Let
$\omega^+, \omega^-\in\Omega$ denote the unique elements representing the
classes $\theta^+$ and
$\theta^-$, respectively.
By Theorem \ref{prelim7}, an element $\eta\in\Omega$ has a decomposition
\begin{equation} \eta=\eta^\perp+c\omega^+
\label{eq7}
\end{equation}
where $\eta^\perp$ belongs to the span of the eigenspaces corresponding to
eigenvalues other than $\rho$.  The fact that $f^*$ and $(f^{-1})^*$ are
adjoint gives
$c=\langle\omega,\theta^-\rangle$.

If $\eta$ is a closed (1,1)-current, then we let $\omega(\eta)$ denote the
element of $\Omega$ that corresponds to the cohomology class $\{\eta\}$
defined by
$\eta$.  Thus $\omega^\pm=\omega(\theta^\pm)$. It is evident that, as a
mapping from currents to
$\Omega$,
$\omega$ is a projection, i.e.,
$\omega\circ\omega=\omega$.   There is a current
$p(\eta)$ such that
\begin{equation}\eta = \omega(\eta)+dd^cp(\eta).
\label{eq8}
\end{equation}
Since $\omega(\eta)$ is
smooth, it follows that $p(\eta)$ is smooth wherever
$\eta$ is.  The potential $p(\eta)$ is uniquely
defined modulo an additive constant, and we specify it uniquely by the
condition
$$\langle p(\eta), \beta\wedge\beta\rangle=0.$$

Now we investigate the interplay between the decomposition (\ref{eq8}) and
$f^*$.  If $\eta$ is positive, then $p(\eta)\in P(X)$, and we may apply
$f^*$ to (\ref{eq8}) to obtain
$f^*\eta=f^*\omega(\eta)+dd^cf^*p(\eta).$
Then we set $$\gamma^+=pf^*\omega$$ and apply the decomposition (\ref{eq8})
to obtain
$$f^*\eta=\omega f^*\omega(\eta)+dd^c[\gamma^+(\eta)+f^*p(\eta)].$$

The operators $\eta\mapsto\omega\eta$ and $\eta\mapsto\gamma^+\eta$ are
linear in
$\eta$ and depend only on the
cohomology class $\{\eta\}$.  The map $\omega$ induces an isomorphism
$\omega:H^{1,1}(X)\to\Omega$.  This provides a conjugacy between the action
of $f^*$ on $H^{1,1}$ and the action of $\omega f^*$ on $\Omega$.  Thus we
have $\omega f^*=\omega f^*\omega$, and we may iterate the previous
equation to obtain
\begin{equation}
f^{n*}\eta=\omega f^{n*}\eta +\rho^ndd^cg_n^+,
\label{eq9}
\end{equation}
where we define
\begin{equation}
g_n^+\eta={1\over \rho^n}\left(f^{*n}p(\eta) +
\sum_{j=0}^{n-1} f^{(n-j-1)*}\gamma^+(f^{j*}\eta)\right).
\label{eq10}
\end{equation}

By Theorem IV.2.13 of [BPV] the set of nef classes, written $H^{1,1}_{\rm
nef}$, has a certain strictness property.  Namely, there is an affine
hyperplane
$H\subset H^{1,1}$
such that
$H\cap H^{1,1}_{\rm nef}$ is compact and generates $H^{1,1}_{\rm nef}$ as a
real cone.  The condition for a class to be nef is
equivalent (see [Lam]) to the condition that its intersections with the
fundamental classes of curves and with the class of the K\"ahler form
$\beta$ are
all non-negative.
If we set
$$
K = \{\eta\in\Omega:\pair{\eta}{V} \geq 0 {\rm\ for\ every\ irreducible\ }
         V\subset\cC(f^{-1})\},
$$
then $K$ is a cone defined by a
finite number of linear inequalities, and it follows from [Lam] that
$H^{1,1}_{\rm nef}\subset K$.

\begin{lem}
\label{prelim29}
Let $K$ be the subset defined above.  Then the function
$$
M(\eta):=\sup_{x\in X-I(f)}
\gamma^+(\eta)x
$$
is finite for $\eta\in K\cap H\cap\Omega$.
\end{lem}

\proof
It is enough to show that for each point $x\in I(f)$ there is a
neighborhood $U$  and a local potential for $f^*\eta$ that is bounded
above.  We have $f^*\eta=\pi_{1*}\pi_2^*\eta$, so if $U\cap I(f)=\{x\}$, we
may argue as in Proposition 1.3 to conclude that on $\pi_1^{-1}U$ we have
$$
\pi_1^* f^*\eta = \pi_2^*\eta + [V]
$$
where $V$ is a (possibly trivial) effective divisor supported on a fiber
$\pi_1^{-1}(x)$.  Hence, $\pi_1^* f^*\eta = dd^c v$ for some
function $v$ on $\pi_1^{-1}(U)$ whose singularities come entirely from local
potentials for $[V]$.  Thus $v$ is bounded above.  Now the pushforward,
$v\circ\pi_1^{-1}$ is a local potential for $f^*\eta$ on $U$ and is
bounded above, as desired.
\qed

\begin{thm}
\label{prelim25}
There are positive constants
$A,B$  such that for any $\eta\in\Omega$
\begin{equation}
\label{prelim20}
|\gamma^+\eta(x)| \leq \norm{\eta}(A+B \left|\log\dist(x,I(f))\right| )
\end{equation}
holds for all $x\in X$.  Further, there exists a constant $C$ such that
if the cohomology class of $\eta$ is nef, then
$\gamma^+\eta \leq C\norm{\eta}$
everywhere on $X$.
\end{thm}

\proof The first assertion follows from writing
$\eta$ as a linear combination of the basis
elements $\omega_1,\dots,\omega_N$ and applying Propositions \ref{prelim1}
and \ref{prelim24}.

By definition, $K$ is a convex cone defined by
finitely many linear inequalities, and since
$H\cap H^{1,1}_{\rm nef}$ is compact, we can choose finitely many
elements
$\eta_1,\dots,\eta_m\in K$ whose convex hull contains $H\cap H^{1,1}_{\rm
nef}$.  The expression
$M(\eta)$ from Lemma
\ref{prelim29} is a convex function of
$\eta$, so we conclude that
$$
C:=\max_{1\leq j\leq m} M(\eta_j)\ge\sup_{\eta\in H\cap H^{1,1}_{\rm
nef}} M(\eta)
$$
gives the upper bound that establishes (\ref{prelim20}).
\qed

\begin{prop}
\label{prelim15}
Given $t>1$, there exists a constant $C$ such that for any
form $\omega\in\Omega$ and any $n\in\N$,
$$
\int |\gamma^+(\omega)|\circ f^j\,dV \leq Ct^j\norm{\omega}.
$$
\end{prop}

\proof From [DF, \S6], we have that for any $t>1$, there exist
constants
$C_1, C_2>0$ such that
$$
\Vol f^{-n}(\ball{I(f)}{r}) \leq C_1 r^{C_2/t^n}
$$
for all $n\in\N$ and all $r>0$.

 From Theorem \ref{prelim25}, we have
$$
\int |\gamma^+\omega|\circ f^j\,dV \leq
\norm{\omega}\left(A + B\int |\log\dist(f^j(x),I(f))|\,dV(x)\right).
$$
Now the volume estimate above gives
\begin{eqnarray*}
\int|\log\dist(f^j(x),I(f))|\,dV(x) & \leq & A + \int_0^\infty
                           \Vol f^{-j}(\ball{I(f)}{e^{-s}})\,ds \\
& \leq &
A + \int_0^\infty C_1 e^{-C_2s/t^j}\,ds
\leq Ct^j,
\end{eqnarray*}
which combines with the first estimate to finish the proof.
\qed

Let us define $\gamma^+:=\gamma^+\omega^+=\gamma^+\theta^+$ so we have
$$dd^c\gamma^+=f^*\omega^+-\rho\omega^+.$$
Since this form is essentially positive, it follows that $\gamma^+$ is
upper semicontinuous.  Thus the infinite sum
\begin{equation}
g^+:=\sum_{j=0}^\infty {\gamma^+\circ f^j\over \rho^j}
\label{eq12}
\end{equation}
is essentially decreasing and defines an upper semicontinuous function
(which is possibly $-\infty$ at some points).
\begin{thm}
\label{prelim14}
The function $g^+$ in (\ref{eq12}) belongs to $L^1(X)$.  Further, for any
smooth, closed
$(1,1)$ form $\eta$, we have
$$
\lim_{n\to \infty} g_n^+\eta = c\cdot g^+
$$
where $c=\langle\eta,\theta^-\rangle$, and the aergence takes place in
$L^1(X)$.
\end{thm}

\proof  Let us first consider the case $\eta=\omega^+$.  Recall that
$\{f^{n*}\omega^+\}=\rho^n\theta^+$ and that $\gamma^+$ depends only on the
cohomology class.  If we set $\gamma^+:=\gamma^+\omega^+=\gamma^+\theta^+$,
then $\gamma^+f^{j*}\omega^+=\rho^j\gamma^+$.  Further, since
$\omega(\omega^+)=\omega^+$, we have $p(\omega^+)=0$.  Thus
$$g^+_n\omega^+ = \sum_{j=0}^{n-1}\frac{\gamma^+\circ f^j}{\rho^j}.$$
If we take
$1<t<\rho$, then by Proposition \ref{prelim15} we have
$$\int|\gamma^+\circ f^j|dV\le C't^{j}.$$
Thus the sequence $\{g_n^+\}$ converges in $L^1(X)$ to
$$g^+:=\sum_{j=0}^\infty \frac{\gamma^+\circ f^j}{\rho^j}.$$

Since $\gamma^+\eta$ depends only on the cohomology class $\{\eta\}$, we
may assume $\eta\in\Omega$.  We use the decomposition  (\ref{eq7}):
$\eta=c\omega^++\eta^\perp$.  Thus
$$g^+_n\eta^\perp={\rho^{-n}}{f^{n*}p\eta^\perp} +
\rho^{-n}\sum_{j=0}^{n-1} f^{(n-j-1)*}\gamma^+(f^{j*}\eta^\perp).$$
Since $\eta^\perp$ is smooth, so is $p\eta^\perp$, and so we have
$|f^{n*}p\eta^\perp|\le C$ on $X$.  Thus $\rho^{-n}f^{n*}p\eta^\perp$
converges uniformly to zero.  By Proposition
\ref{prelim15} again and the fact that $\gamma^+=\gamma^+\circ\omega$, we
have
$$\int|g^+_n\eta^\perp|dV\le C\frac{{\rm Vol}(X)}{\rho^n} +
\frac1{\rho^n}\sum_{j=0}^{n-1} Ct^{n-j-1}||\omega f^{j*}\eta^\perp||.$$
By Theorem \ref{prelim7}, there is a constant $C'$
such that
$||\omega f^{j*}\eta^\perp||\le C't^j$.   Thus
$$\int|g^+_n\eta^\perp|dV\le C\frac{{\rm Vol}(X)}{\rho^n} +
  CC'n\frac{t^{n-1}}{\rho^n}.$$
This tends to zero as $n\to\infty$, and $g^+_n\eta$ is linear in $\eta$,
so the Theorem follows.  \qed

\begin{thm}
\label{prelim17}
The current $\mu^+ := \omega^+ + dd^c g^+$ has the following properties:
\begin{itemize}
\item for every smooth closed $(1,1)$ form $\eta$ on $X$, we have
$$
\lim_{n\to\infty} \frac{f^{n*}\eta}{\rho^n}
            = \mu^+\cdot\pair{\eta}{\theta^-}.
$$
\item $\mu^+$ is positive;
\item $f^*\mu^+ = \rho\mu^+$.
\end{itemize}
\end{thm}

\proof
By Theorem \ref{prelim14},
$$
\lim_{n\to\infty}\frac{f^{n*}\eta}{\rho^n}
            = \lim_{n\to\infty} \frac{\omega f^{n*}\eta}{\rho^n}
              + dd^c \lim_{n\to\infty} g^+_n\eta
            = \pair{\eta}{\theta^-}(\omega^+ + dd^c g^+).
$$
Taking $\eta = \beta$ to be the K\"ahler form on $X$, we have
$\pair{\beta}{\theta^-} = 1$.  Thus
$
\mu^+
       = \lim_{n\to\infty}{\rho^{-n}}{f^{n*}\beta}
$
is a limit of positive currents and therefore positive.
Since $f^*$ acts continuously on positive closed $(1,1)$ currents, we
also get that
$
f^*\mu^+ = \lim_{n\to\infty}{\rho^{-n}}{f^{(n+1)*}\beta} = \rho\mu^+.
$
\qed

The following is an observation of Favre [Fav1, Theorem 1]:
\begin{cor}
\label{prelim228}
The Lelong number $\nu(\mu^+,x)$ vanishes for
$x\in X-\bigcup_{n\ge1}I(f^n)$.
\end{cor}

Now we discuss the extent to which $\mu^+$
is invariant under bimeromorphic conjugacy.
\begin{prop}
Let $\pi:\tilde X\to X$ be a proper modification of $X$.  Suppose that
$\tilde f=\pi^{-1}f\pi$ is a bimeromorphic map of $\tilde X$ which
satisfies (\ref{eq1}) and (\ref{eq5}).  If $\mu^+$ and $\tilde\mu^+$
denote the associated currents, then $\pi_*\tilde\mu^+$ is a positive
multiple of
$\mu^+$.
\end{prop}

\proof
By hypothesis $I(\pi) = \emptyset = \cC(\pi^{-1}) = \emptyset$.
Thus $\pi^*\beta$ is a smooth, positive and closed $(1,1)$ form on $X$,
and we compute
$$
c' \pi_* \tilde\mu^+
    = \pi_*\lim_{n\to\infty} \frac{\tilde f^{n*} (\pi^*\beta)}{\rho^n} \\
    = \lim_{n\to\infty} \frac{\pi_* \tilde f^{n*} \pi^*\beta}{\rho^n} \\
   = \lim_{n\to\infty} \frac{ f^{n*}\beta }{\rho^n} = \mu^+
$$
The first and last equalities follow from Theorem
\ref{prelim17}.  The second inequality follows from continuity of $\pi_*$
acting on positive closed $(1,1)$ currents, and the third equality is
a consequence of the proof of Proposition 1.13 in [DF].
Clearly $c'>0$.
\qed

\begin{prop}
\label{prelim6}
Let $h:X\to Y$ be a bimeromorphic map, and let
$\tilde f=h fh^{-1}$ be the induced bimeromorphic self-map of
$Y$.  If
$f$ and
$\tilde f$ satisfy (\ref{eq1}) and (\ref{eq5}), then
$h^*\tilde\mu^+ = ch^*\mu^+ + [V]$,
where $c>0$ and $V$ is an effective divisor supported on $\cC(h)$.
\end{prop}

\proof
Let $G$ be the desingularized graph of $h$ and
$\pi_X:G\to X$, $\pi_Y:G\to Y$ be the projections onto first and second
factors.  That is, $h = \pi_Y\circ \pi_X^{-1}$.  After blowing up
points in $G$ if necessary, we can assume that the common lift
$F:G\to G$ of $f$ and $\tilde f$ to $G$ satisfies (\ref{eq1}).
Let $\nu^+$ denote the expanding current associated to $F$.  Then
by the previous lemma, we see that
$$
\pi_{X*}\nu^+ = c_1\mu^+ \qquad \pi_{Y*}\nu^+ = c_2\tilde\mu^+
$$
for constants $c_1, c_2 > 0$.  Hence, $c\pi_X^*\mu^+ - [V'] =
\pi_Y^*\tilde\mu^+ - [V'']$ where $c>0$ and $V'$ and $V''$ are effective
divisors supported on $\cC(\pi_X)$ and $\cC(\pi_Y)$, respectively.  We
apply the `pushpull formula' [DF, Theorem 3.3] to conclude
$$
c\mu^+  = c\pi_{X*}\pi^*_X \mu^+ = c\pi_{X*}(\pi_X^*\mu^+ - [V']) =
h^*\tilde\mu^+ - [\pi_{X*}V''].
$$
Since $\pi_{X*}V''$ is an effective divisor supported on $\cC(h)$,
we are done.
\qed

\section{Energy}
\label{energy}
Let $T$ be a positive, closed $(1,1)$ current on $X$.  Then $T$ defines
an inner product on the space of smooth, real functions on $X$ via the
formula
$$
\cE_T(\varphi,\psi) \eqdef \int d\varphi\wedge d^c\psi\wedge T.
$$
We denote the seminorm associated with this inner product by
$$|\varphi|_T=\cE(\varphi,\varphi)^{1/2}=\left(\int d\varphi\wedge
d^c\varphi\wedge T\right)^{1/2}.$$

We will say that functions $u_j\in C^\infty(X)$, $j\ge0$ form
a \emph{regularizing sequence} for a function $u$ if $u_j$ decreases pointwise
to $u$ and $dd^c u_j \geq - c\beta$ for some $c>0$ and all $j$.  The
limit $u$ necessarily belongs to $P$, and indeed
any $u\in P$ admits a regularizing
sequence (see [Dem1, Theorem 1.1]).  We will use the following property of
a function $u\in P$:
\begin{equation}
{\rm\ Every\ regularizing\ sequence\ }\{u_j\}{\rm\ for\ }u{\rm\ is\ Cauchy\
in\ }|\cdot|_T.\label{tregular}
\end{equation}
The union of two regularizing sequences is (essentially) a regularizing
sequence.  Thus if $u$ satisfies (\ref{tregular}) then all regularizing
sequences define the same element of the completion with respect to
$|\cdot|_T$.  In particular, if
$u$ and
$v$ satisfy (\ref{tregular}), then we may define
$\cE_T(u,v)$ by taking the limit along regularizing sequences.

The special case $T=\beta$ is classical: condition (\ref{tregular})
for $T=\beta$ is equivalent to the condition that $\nabla u\in L^2$.

Our principal use of
condition (\ref{tregular}) is to define $(dd^cu)\wedge T$.  If $u$
satisfies (\ref{tregular}), then we may define $(dd^cu)\wedge T$ as a
distribution via the pairing
$$\psi\mapsto\langle dd^cu\wedge T,\psi\rangle:= -\cE_T(\psi,u).$$
It is evident that $dd^cu\wedge T+c\beta\wedge T\ge0$, so $dd^cu\wedge T$ is
represented by a (signed) Borel measure.  Further, since $dd^cu_j\ge-c\beta$,
it also follows from (\ref{tregular}) that $dd^cu_j\wedge T$ converges to
$dd^cu\wedge T$ in the weak* topology on the space of measures.

\begin{prop}
If $u\in L^1(T\wedge\beta)$, and if $u$ satisfies (\ref{tregular}), then
$(dd^cu)\wedge T=dd^c(uT)$.
\end{prop}

\proof  Let us remark first that if $u\in L^1(T\wedge\beta)$, then $uT$ is a
well-defined (1,1) current, and thus $dd^c(uT)$ is a well-defined current.  If
$\{u_j\}$ is a sequence satisfying (\ref{tregular}), then $u_jT$ converges to
$uT$ weakly as currents.  Thus $dd^c(u_jT)$ converges to $dd^c(uT)$.  Finally,
$(dd^cu_j)\wedge T=dd^c(u_jT)$ when $u_j$ is smooth, and we have observed
above that $\lim_{j\to\infty}(dd^cu_j)\wedge T=(dd^cu)\wedge T$. \qed

\begin{prop}
\label{prop32}
If $u,v\in P$ both satisfy (\ref{tregular}), and if $v\in L^1(T\wedge\beta)$,
then
$v\in L^1(dd^cu\wedge T)$.
\end{prop}

\proof Let $\{u_j\}$ and $\{v_k\}$ denote regularizing sequences for $u$ and
$v$.  For fixed $j$ and $k$, integration by parts gives $\int
v_k\,dd^cu_j\wedge T=-\int dv_k\wedge d^cu_j\wedge T$.  Now
$dd^cu_j+c\beta\ge0$, so $(dd^cu+c\beta)\wedge T$ defines  (positive) Borel
measure.  Letting
$j\to\infty$, we have
$$\int|v_k|\, (dd^cu+c\beta)\wedge T = -\cE_T(v_k,u) + c\int|v_k|\,\beta\wedge
T.$$
If we let $k\to\infty$, then the right hand side stays bounded since $v\in
L^1(T\wedge\beta)$, and thus $v\in L^1(dd^cu\wedge T)$ by monotone
convergence.
\qed

The motivation for our work in the following sections is as follows.  We will
show that $g^+\in L^1(T\wedge\beta)$, and
$g^+$ satisfies (\ref{tregular}) for the current
$T=\mu^-+\beta$.  It will then follow that
$dd^cg^+\wedge\mu^-$ is well defined, so the wedge product defines a
(signed) measure
\begin{equation}\mu=\mu^+\wedge\mu^-=\omega^+\wedge\mu^-+dd^cg^+\wedge\mu^-,
\label{defofmu}
\end{equation}
and $g^+\in L^1(\mu)$.  Since $\mu^\pm\ge0$, it follows that $\mu$ is positive.
The total mass of $\mu$ is 
$\int\omega^+\wedge\omega^-=\langle\theta^+,\theta^-\rangle=1$,
so $\mu$ is a probability measure.

\begin{lem}
\label{nrg8}
Let $u,v\in C^\infty(X)$ satisfy $dd^c u,dd^c v\geq -c\beta$ and
$v\geq u$.  Then for any positive, closed $(1,1)$ current $T$,
\begin{eqnarray*}
\cE_T(u,v) - \cE_T(v,v) & \geq & -c \int(v-u)\,\beta\wedge T \\
\cE_T(u,u) - \cE_T(u, v) & \geq & -c \int(v-u)\,\beta\wedge T.
\end{eqnarray*}
\end{lem}

\proof
It is sufficient to prove the first inequality.
\begin{eqnarray*}
\cE_T(u,v)-\cE_T(v,v)
& = &
\int d(u-v)\wedge d^cv\wedge T\\
& = &
\int(v-u)dd^cv\wedge T
\geq
-c\int(v-u)\beta\wedge T\\
\end{eqnarray*}
Here we used Stokes' Theorem to pass from the first line to the second
line, and the inequality is obtained because $v-u\ge0$ and
$dd^cv\geq -c\beta$.
\qed

\begin{thm}
\label{nrg7}
Suppose that $T$ is a positive closed $(1,1)$ current and
$u\in L^1(\beta\wedge T)$.  If there exists a regularizing sequence $\{u_j\}$
for which $\{|u_j|_T\}$ is bounded, then
$u$ satisfies (\ref{tregular}).
\end{thm}

\proof
By hypothesis there exists $c>0$ such that
$dd^cu_j+c\beta\geq 0$ for all $j$.  Now $\beta\wedge T$ is a positive,
finite Borel measure, and by the monotonicity of the sequence $\{u_j\}$, we
have
$$
\lim_{j\to\infty}\int|u_j-u|\,\beta\wedge
T = \lim_{j,k\to\infty}\int|u_j-u_k|\,\beta\wedge T=0.
$$
It follows from
the Lemma \ref{nrg8} that for $k\ge j$, we have
\begin{eqnarray*}
\cE_T(u_k,u_k)&-&\cE_T(u_j,u_j)\\
&\ge&
\cE_T(u_k,u_k)-\cE_T(u_k,u_j) + \cE_T(u_k,u_j)-\cE_T(u_j,u_j)\\
&\ge&
-2c\int|u_k-u_j|\,\beta\wedge T
\end{eqnarray*}
Thus the sequence $|u_k|_T=\cE_T(u_k,u_k)^{1/2}$ is essentially
increasing.  Since we have assumed that it is also bounded, we
conclude $\lim_{k\to\infty}|u_k|_T$ exists and is finite.

Now we observe that
$$
|u_j-u_k|_T^2
=
\cE_T(u_k,u_k)-2\cE_T(u_j,u_k)+\cE_T(u_j,u_j)
$$
If $k\ge j$, then $u_k\le u_j$, so by the Lemma, we have $\cE_T(u_j,u_j)\le
\cE_T(u_j,u_k)\le
\cE_T(u_k,u_k)$, modulo an error of size $2c\int|u_k-u_j|\omega\wedge T$.
Thus
$$\lim_{j\to\infty}\cE_T(u_j,u_j)=\lim_{j,k\to\infty}\cE_T(u_j,u_k)=
\lim_{k\to\infty}\cE_T(u_k,u_k)$$
and so $\lim_{j,k\to\infty}|u_j-u_k|_T=0$.

Now we show that (\ref{tregular}) holds.  Let $\{v_j\}$ be any regularizing
sequence for $u$.  Since $v_j$ is smooth, there exists $k=k_j$ such that
$u_{k_j}\le v_j$.  Thus by Lemma \ref{nrg8}, $|v_j|_T$ is essentially bounded
by $|u_{k_j}|_T$.  From the first part of the proof, then, it follows that
$\{v_j\}$ is Cauchy.
\qed

\begin{prop}
\label{nrgprop}
Suppose that $u\in\tilde P(X,S)$ and that $T$ is a positive closed
$(1,1)$ current whose local potential is finite at each point of $S$.
Then $u\in L^1(T\wedge\beta)$.
\end{prop}

\proof
First recall from Jensen's inequality that
$[\beta\wedge T](\ball{S}{r}) \leq Cr^2$ for all $r\geq 0$.  Since
$u(x) \geq A\log\dist(x,S) + B$, we have that
$$
[\beta\wedge T] \{|u|\geq t\} \leq [\beta\wedge T](\ball{S}{e^{(-t+B)/A})}
\leq C'e^{-2t/A}
$$
for constants $C,C'>0$ and all $t\geq 0$.  Therefore,
$$
\int |u|\,\beta\wedge T = \int_0^\infty [\beta\wedge T]\{|u|\geq t\}\,dt
                         \leq C'\int_0^\infty r^{-2t/A}\,dt < \infty.
$$
\qed

\begin{thm}
\label{nrg5}
Suppose that $u\in\tilde P(X,S)$ and that $T$ is a positive closed
$(1,1)$ current with local potentials that are finite at each point in $S$.
Then $u$ satisfies (\ref{tregular}), and so $|u|_T:=\lim_{j\to\infty}|u_j|_T$
is well defined and finite; and
$$|u|_T^2=\int_{X-S}du\wedge d^cu\wedge T.$$
The  expression
$du\wedge T$ defines a current on
$X-S$ which has finite mass.  The (trivial) extension of $du\wedge T$ to $X$
is equal to the current $d(uT)$, i.e., for all smooth 1-forms $\eta$ we have
$$
\pair{d(u T)}{\eta} = \int_{X\setminus S} \eta\wedge du\wedge T.
$$
\end{thm}

\proof
First we show that $u$ satisfies (\ref{tregular}).  Choose a
function $m\in C^\infty(\R)$
that is convex, increasing and equal to $\max\{0,t\}$ outside a small
neighborhood of $0$.  For all $j\ge0$ let
$u_j(x) = m(u(x) + j) - j\approx \max\{u(x),-j\}$.  Clearly $u_j$
is smooth and decreases to $u$ pointwise on $X$.  An easy computation
verifies that
if $dd^c u\geq - c\beta$, then the same is true for $u_j$ with the
same constant $c$.  That is, $u_j$ regularizes $u$.  By Theorem \ref{nrg7}
and Proposition \ref{nrgprop} it suffices to show that $|u_j|_T$ is bounded.

Let $L\in C^\infty(X\setminus S)$ be a function satisfying $L\leq u$ and
$L(z) = C(q)\log\norm{z}$ with respect to local coordinates centered at
each point $q\in S$.  For each $j\in\N$, let $L_j(z) =  m(L(z) + j) - j$.
Then the above argument applied to $L,L_j$ instead of $u,u_j$ shows
that $L\in L^1(\beta\wedge T)$ and that $L_j$ regularizes $L$. Since
$u_j\geq L_j$, we have from Lemma \ref{nrg8} that
$$
|u_j|_T\leq |L_j|_T + C \int (u_j - L_j)\,\beta\wedge T
\leq |L_j|_T + C\int (\norm[\infty]{u} - L)\,\beta\wedge T
\leq |L_j|_T + C
$$
for every $j\in\N$.  Hence our problem reduces to showing that
$\{|L_j|_T\}$ is bounded.

To do this, we can restrict attention to a coordinate neighborhood
$\ball{0}{2}$ centered at $q\in X$ and assume that
$L(z) = C\log\norm{z}$ in these coordinates.  We choose a smooth,
radially symmetric and compactly supported function
$\chi:\ball{0}{2}\to [0,1]$ such that $\chi\equiv 1$ on $\ball{0}{1}$.
We choose a local potential $v$ for $T$ on $\ball{0}{2}$ and estimate
$$
\int_{\ball{0}{1}} dL_j(z)\wedge d^cL_j(z)\wedge T
\leq
\int \chi\,dL_j(z)\wedge d^cL_j(z)\wedge dd^c v.
$$
Integrating by parts, we see that the right side is dominated by
\begin{eqnarray*}
\left|\int v \,dd^c\chi\wedge dL_j\wedge d^cL_j\right| &+&
          2\left|\int v \,d\chi\wedge d^cL_j\wedge dd^c L_j\right| \,\,+\,\,
          \left|\int \chi v\,dd^cL_j\wedge dd^c L_j\right| \\
  = \left|\int v \,dd^c\chi\wedge dL\wedge d^cL\right| &+&
          2\left|\int v \,d\chi\wedge d^cL\wedge dd^c L\right| \,\,+\,\,
          \left|\int \chi v\,dd^cL_j\wedge dd^c L_j\right|
\\
& = & C + \int (-v)\chi \,dd^cL_j\wedge dd^c L_j
\end{eqnarray*}
for $j$ large enough that $L_j = L$ on $\supp d\chi$.  The
measures $\chi(dd^c L_j)^2$ are radially symmetric and converging to
a point mass at the origin as $j\to\infty$.  Since $v$ is subharmonic,
we obtain that
$$
- \lim_{j\to\infty} \int \chi v\,dd^cL_j\wedge dd^c L_j = -v(0) <\infty.
$$
Thus $\{|L_j|_T\}$ is bounded, and (\ref{tregular}) holds by Theorem
\ref{nrg7}.

The Radon-Nikodym Theorem allows us to write  
$$
du_j\wedge d^cu_j\wedge T = g_j \beta\wedge T,\qquad 
(du\wedge d^cu\wedge T)|_{X-S} = g\beta\wedge T,
$$
where $g,g_j\in L^1(\beta\wedge T)$, and by the above work, $\{g_j\}$ 
converges in $L^1(\beta\wedge T)$ to a function $g_\infty$ that is 
independent of the regularizing sequence $u_j$.
Since we can choose $\{u_j\}$ to be equal to $u$ outside of an arbitrarily
small neighborhood of $S$, it follows that $g_\infty = g$ outside $S$.  But
$T\wedge \beta$ does not charge $S$, so the formula for $|u|_T$
follows.

To prove the formula for $\pair{d(uT)}{\eta}$, let $\chi,1-\chi$ be a partition
of unity subordinate to $U,X-K$, where $U$ is a small neighborhood of $S$ and
$K\subset U$ is a closed set containing a neighborhood of $S$.  Then
$$
\pair{d(uT)}{\eta} = \pair{d(uT)}{\chi\eta} + 
	\int (1-\chi) \eta\wedge du \wedge T.
$$
Since $u_jT\to uT$ weakly, we can apply the Schwarz inequality to obtain
$$
|\pair{d(uT)}{\chi\eta}| = \lim_{j\to\infty} |\pair{d(u_jT)}{\chi\eta}|
	\leq \lim_{j\to\infty} |u_j|_T 
	\left(\int \chi^2\eta\wedge J\eta\wedge T\right)^{1/2}
	\leq C|u|_T \int_{U} \beta\wedge T	
$$
which goes to zero as the neighborhood $U$ shrinks to $S$.  The formula
for $\pair{d(uT)}{\eta}$ follows.
\ignore{Finally, in order to prove the first equation in the statement of the
Theorem it suffices to show that for $x\in S$ we may choose $r$ small so
that the integral
$\int_{B_x(r)}du_j\wedge d^cu_j\wedge T$ is arbitrarily small, uniformly
in
$j$.  Arguing as above, we may assume that $x=0$ and replace $u_j$ by the
logarithm
$L(z)$, so it suffices to show that $\int_{B_x(r)}du_j\wedge d^cu_j\wedge
T$ is uniformly small in $j$.  In local coordinates near $x=0$, we have
$dL\wedge d^cL\le c||z||^{-2}\beta$.   We may identify the
Laplacian of the local potential $v$ for $T$ as the measure
$\Delta v = T\wedge \beta$.  Since
$-||z||^{-2}$ is the Newtonian potential on ${\bf R}^4$,  we may assume that
$v$ is given by convolution: $v(0)=-v\star ||z||^{-2}(0)=-\int||z||^{-2}\Delta
v$. Since
$v(0)$ is finite, it follows that $||z||^{-2}$ is integrable with respect to
the measure $\Delta v$.  Thus the integral over $B_x(r)$ may be taken
uniformly small in $j$.}
\qed

\begin{cor}
\label{nrg6}
Let $u\in\tilde P(X,S)$ where $S\cap I(f^{-1}) = \emptyset$, and $T$ be a
positive closed $(1,1)$ current with
local potentials that are finite at each point in $I(f)\cup f^{-1}(S)$.  Then
$|f^*u|_T$ is well-defined, and
\begin{itemize}
\item $\int_{X-S}du\wedge d^cu\wedge f_*T<\infty$
\item $|f^*u|^2_T=\int_{X-(I(f)\cup f^{-1}S)}d(f^*u)\wedge d^c(f^*u)\wedge
T<\infty$
\item $|u|_{f_*T} = |f^*u|_T.$
\end{itemize}
\end{cor}

\proof By Proposition \ref{prelim24}, $f^*u=v_1-v_2$ is the difference of
functions in $\tilde P(X,I(f)\cup f^{-1}S)$.  By Theorem \ref{nrg5}, we may
define $|v_j|_T$ in terms of the integral of the pointwise gradient $dv_j$
over $X-(I(f)\cup f^{-1}S)$.  Thus $|f^*u|_T$ is well defined. The hypotheses
imply that
$f_*T$ does not charge
$\cC(f^{-1}) = f(I(f))$ and  that
$f_* T$ has local potentials that are finite at each point
$p\in S \subset f(f^{-1}(S))$.  The function $u\circ f$ is a difference
of elements of $\tilde P(X,f^{-1}(S)\cup I(f))$ by
Proposition \ref{prelim24}.   Hence by Theorem
\ref{nrg5}, the integrals defining $|u|_{f_*T}$ and $|u\circ f|_T$ are
finite.   We compute
\begin{eqnarray*}
|u|_{f_*T}^2
& = &
\int_{X - I(f^{-1}) - \cC(f^{-1})} du\wedge d^c u\wedge f_*T \\
& = &
\int_{X - I(f) - \cC(f)} d(u\circ f)
                  \wedge d^c(u\circ f)\wedge T \\ & = &
\int_{X - I(f)} d(u\circ f)\wedge d^c(u\circ f)\wedge T =
|u\circ f|_T^2.
\end{eqnarray*}
The first equality holds because $f_*T$ charges neither points nor
$\cC(f^{-1})$.
The second equality follows by the change of variables formula because
$f:X-I(f)-\cC(f) \to X-I(f^{-1})-\cC(f^{-1})$ is a biholomorphism.
The third equality is a consequence of the fact that $u\circ f$ is constant
on $\cC(f)- I(f)$.  Finally, the fourth equality holds because
$T$ does not charge points.
\qed

\section{Invariant Measure}
\label{mu}
Up to this point, we have required that conditions (\ref{eq1}) and 
(\ref{eq5}) hold.  We
will now impose two further conditions.  The first of these is:
\begin{equation}
\langle\theta^+,f(x)\rangle>0{\rm\ for\ every\ }x\in I(f),{\rm\ and\ }
\langle\theta^-,f^{-1}(y)\rangle>0{\rm\ for\ every\ }y\in I(f^{-1}).
\label{eq13}
\end{equation}
Like condition (\ref{eq1}), condition (\ref{eq13}) may be thought of 
as a property of the
underlying space $X$ used to represent the map $f$; it will be shown 
in Proposition
\ref{prelim18} that (\ref{eq13}) may always be assumed to hold. 
Next we consider
condition (\ref{eq3}) more carefully (Theorems \ref{mu15} and 
\ref{mu1}).  After this, we
will assume for the rest of the paper that (\ref{eq3}) holds, by which we
mean implicitly that (\ref{eq1}), (\ref{eq3}), (\ref{eq5}) and 
(\ref{eq13}) all hold.  The
main results of this section are that if (\ref{eq3}) holds, then the 
expression $\mu$ in
(\ref{defofmu}) is well defined (Theorem \ref{mu2}) and invariant (Theorem
\ref{invariance}).

\begin{prop}
\label{prelim18}
If $f:X\to X$ satisfies (\ref{eq1}), we may blow down curves in $X$ if
necessary so that both (\ref{eq1}) and (\ref{eq13}) hold.
\end{prop}

\proof
Suppose to the contrary that $\pair{\theta^+}{f(x)} \le 0$ for some $x\in I$.
Then since $\theta^+$ is nef, we have $\pair{\theta^+}{V} = 0$ for every
component $V\subset f(x)$.  From this and the Hodge index theorem on
surfaces, we see that either $\pair{\theta^+}{\theta^+} = 0$ or the
intersection
form is negative definite on $f(x)$.  In the first case
[DF, Theorem 0.4]
guarantees that after blowing down an appropriate curve, $f$ conjugates to an
automorphism and satisfies the conclusion of the proposition vacuously.

In the second case, we note that $\pi^{-1}_1(x)$ constitutes a single
connected component of $\cC(\pi_1)$.  We can therefore apply the argument
of [DF, Proposition 1.7] to obtain a smooth rational curve of
self-intersection
$-1$ in $f(x)$.  After blowing this curve down, (\ref{eq1}) still
holds.  However the dimension of $H^{1,1}(X)$ drops by one.  If on the
new surface we still have $\pair{\theta^\pm}{f(x)} = 0$ for some 
$x\in I(f^{\pm1})$,
then we can repeat this process.  This cannot happen more than
$\dim H^{1,1,}(X)$ times, so eventually we will descend to a
surface on which $f$ satisfies both (\ref{eq1}) and (\ref{eq13}).
\qed

The following is a companion to Corollary \ref{prelim228}.
\begin{cor}
\label{lelongnumber} If $f$ satisfies (\ref{eq13}), then $\nu(\mu^+,x)
>  0$ for $x\in I(f^n)$, $n\ge1$.
\end{cor}

\proof The proof uses the fact that
$\mu^+ = \rho^{-n}f^{n*}\mu^+$, but it is otherwise identical to the
proof of the second conclusion in Proposition \ref{prelim1}.
\qed

\begin{thm}
\label{mu15}
Suppose that (\ref{eq1}), (\ref{eq5}), and (\ref{eq13}) hold.  Then 
(\ref{eq3}) holds
if and only if the function $g^+$ defined in (\ref{eq12}) is finite 
at each point of
$I(f^{-1})$.
\end{thm}

\proof
Remarking that $f^*\omega^+=\rho\omega^+ +dd^c\gamma^+$, we see from 
Proposition 1.3 and
Corollary \ref{lelongnumber} that
$$A\log \dist(x,I(f))-B\le\gamma^+(x)\le A'\log \dist(x,I(f))+B'.$$
Replacing $x$ by $f^jx$ and summing we see that $g^+$ is bounded 
above and below by
infinite sums of the form  $S:=A\sum \rho^{-n}\log \dist(f^nx,I(f))+B.$   Thus
$g^+(x)>-\infty$ if and only if $S$ is finite.  Since
$I(f^{-1})$ is a finite set, we have
$$\dist(f^nI(f^{-1}),I(f))=\min_{x\in I(f^{-1})} \dist(f^jx,I(f)),$$
and we see that (\ref{eq3}) holds if and only if $S$ is finite for all $x\in
I(f^{-1})$. \qed

The following Theorem says that condition (\ref{eq3}) is symmetric in $f$ and
$f^{-1}$; we refer the reader to [Dil1, Theorem 5.2] for a proof.
\begin{thm}
\label{symmetry}
Suppose that (\ref{eq1}), (\ref{eq5}), and (\ref{eq13}) hold.  Then 
(\ref{eq3}) holds
if and only if
\begin{equation}
\label{mu16}
\sum_{j=0}^\infty \frac{\log\dist(f^{-j}(I(f)),I(f^{-1}))}{\rho^j} < \infty.
\end{equation}
\end{thm}

\begin{prop}
\label{mu1}  If (\ref{eq3}) holds, then there exists a constant $C>0$ 
such that for
each $\eta\in\Omega$ and every $n\in\N$
$$
|f^{n*}\gamma^+\eta|_{\beta+\mu^-} \leq C\rho^{n/2}\norm{\eta}.
$$
\end{prop}

\proof
The function $\gamma^+\eta\circ f^n$ is a difference of functions in
$\tilde P(X,I(f^n))$ by Theorem \ref{prelim25}, and $f^n_*\beta$
is smooth away from $I(f^{-n})$, so by Corollary \ref{nrg6}
$f^{n*}\gamma^+\eta$ is a difference of functions which satisfy condition
(\ref{tregular}) for $T={f^m_*\beta}$.  The condition (\ref{eq3}) together with
Proposition \ref{mu15} imply that
$f^{n*}\gamma^+\eta$ is also a difference of functions satisfying 
(\ref{tregular}) for
$T=\mu^-$.

 From Corollary \ref{nrg6} again and the invariance of $\mu^-$, we have
$$
|f^{n*}\gamma^+\eta|_{\mu^-} =
|\gamma^+\eta|_{f^n_*\mu^-} =
\rho^{n/2}|\gamma^+\eta|_{\mu^-}\leq C\rho^{n/2}\norm{\eta},
$$
where the last inequality follows from writing $\eta$ as a linear
combination of the basis elements $\omega_1,\dots,\omega_n\in\Omega$.

Corollary \ref{nrg6} also gives
$$|f^{n*}\gamma^+\eta|_\beta 
=\rho^{n/2}|\gamma^+\eta|_{\rho^{-n}f^n_*\beta} \le
\rho^{n/2}||\eta||\max_{k}|\gamma^+\omega_k|_{\rho^{-n}f^n_*\beta}.$$
It will suffice to show that 
$|\gamma^+\omega_k|_{\rho^{-n}f^n_*\beta}$ is bounded.
 From \S2, we write
$$\rho^{-n}f^n_*\beta = \rho^{-n}\omega f^n_*\beta + dd^cg^-_n\beta.$$
With this notation, we have
\begin{equation}|\gamma^+\omega_k|^2_{\rho^{-n}f^n_*\beta} =
|\gamma^+\omega_k|^2_{\rho^{-n}\omega f^n_*\beta} +
|\gamma^+\omega_k|^2_{dd^cg^-_n\beta}.\label{eq13a}
\end{equation}
The sequence $\{\rho^{-n}\omega f^n_*\beta: 
n=1,2,3,\dots\}\subset\Omega$ is bounded and in
fact converges to
$c\omega^-$.  Thus the first term on the right hand side is bounded.

Before we analyze the
second term, we make some observations.  Since $\gamma^+\omega_k\in 
P$ is smooth except for
logarithmic singularities at
$I(f)$, the wedge product yields a well-defined measure
$$(dd^c\gamma^+\omega_k)^2=\sum_{x\in I(f)} c_x\delta_x +
(dd^c\gamma^+\omega_k)^2|_{X-I(f)}$$
(see [Dem2]).  We claim that the integral $\int
(g_n^-\beta)(dd^c\gamma^+\omega_k)^2$ is uniformly bounded in $n$ if 
(\ref{eq3}) holds.
To see this we use the projections
$\pi_1,\pi_2:\Gamma\to X$ and the formula 
$f^*\omega_k=\pi_{1*}\pi_2^*\omega_k$.  Since
$\pi_1:\Gamma-\cC(\pi_1)\to X-I(f)$ is biholomorphic, we may pull 
back by $\pi_1^*$ to
obtain:
$$\int_{X-I(f)} (g^-_n\beta)(dd^c\gamma^+\omega_k)^2 = \int_{X-I(f)}
(g^-_n\beta)(f^*\omega-\rho\omega_k)^2 = \int_\Gamma
\pi^*_1(g^-_n\beta)(\pi_2^*\omega_k-\rho\pi_1^*\omega_k)^2.$$
Now $(\pi_2^*\omega_k-\rho\pi_1^*\omega_k)$ is a smooth (1,1)-form on 
$\Gamma$, and
$\pi_1^*(g^-_n\beta)$ is a sequence which decreases monotonically to 
$\pi_1^*g^-_n\beta\in
L^1(\Gamma)$.  Thus this family of integrals is uniformly bounded in 
$n$.  The other
integral,
$$\int_{I(f)}(g^-_n\beta)(dd^c\gamma^+\omega_k)^2 =\sum_{x\in 
I(f)}c_x(g^-_n\beta)(x)$$
is uniformly bounded by (\ref{eq3}) and Proposition 4.1.

With these observations we will now show that the second term in 
(\ref{eq13a}) is bounded.
First we claim that
$$|\gamma^+\omega_k|^2_{dd^cg^-_n\beta}= \int d(\gamma^+\omega_k)\wedge
d^c(\gamma^+\omega_k)\wedge dd^cg^-_n\beta = -\int 
(g^-_n\beta)(dd^c\gamma^+\omega_k)^2.$$
To see this, replace $\gamma^+\omega_k$ by a regularizing
sequence $\gamma_j$, which coincides with $\gamma^+\omega_k$ outside 
neighborhood of
$I(f)$.  Performing two integrations by parts, we have
$$|\gamma_j|_{dd^cg^-_n\beta}=-\int (g^-_n\beta)(dd^c\gamma_j)^2.$$
Now we take the limit as
$j\to\infty$.  The measures $(dd^c\gamma_j)^2$ converge weakly to
$(dd^c\gamma^+\omega_k)^2$.  Further, since $g^-_n$ is continuous in 
a neighborhood of
$I(f)$, we may assume that it is continuous on the set where
$\gamma_j\ne\gamma^+\omega_k$.  This establishes the claim.
Letting $n\to\infty$ and appealing to condition (\ref{eq3}), we see 
that the second term
on the right hand side of (\ref{eq13a}) is bounded.
\qed

Now we sharpen Theorem \ref{prelim14}.

\begin{thm}
\label{mu2}
If (\ref{eq3}) holds, then $g^+\in L^1(\mu^-\wedge\beta)$, and $g^+$ satisfies
(\ref{tregular}) for $T={\beta+\mu^-}$.  In fact, if
$\eta$ is any smooth closed (1,1) form on $X$, then
$$
\lim_{n\to\infty}|g^+_n\eta-g^+|_{\beta+\mu^-}=0.
$$
\end{thm}

\proof  First we will show that for every smooth, closed (1,1) form 
$\eta$, the sequence
$\{g_n^+\eta\}$ is Cauchy.  We consider increasingly general forms 
$\eta$.  First, if
$\eta=\omega^+$, it follows from Proposition 4.3 that 
$\{g^+_n\omega^+\}$ is Cauchy in
the seminorm
$|\cdot|_{\beta+\mu^-}$.  Now for  $\eta\in\Omega$, we write
$\eta=c\omega^++\eta^\perp$ as in (\ref{eq7}).  By Theorem 
\ref{prelim7}, we have
$||\omega f^{j*}\eta^\perp||\le C_tt^j||\eta^\perp||$ for any $t$ 
greater than 1.
By Proposition 4.3 again, we have
$$|\rho^{-n}f^{(n-j-1)*}\gamma^+f^{j*}\eta^\perp|_{\beta+\mu^-}\le
C\rho^{-n}\rho^{(n-j-1)/2}t^{j}.$$ Choosing $t<\sqrt\rho$, we see that
$|g_n^+\eta^\perp|_{\beta+\mu^-}$ converges to zero.  Since 
$g^+_n\eta$ is linear in
$\eta$, we see that $g^+_n\eta$ satisfies (\ref{tregular}) with 
$T=\beta+\mu^-$.
Finally, for a general, smooth, closed (1,1) form $\eta$, we see that 
$g_n^+\eta$ differs
from
$g_n^+\omega\eta$ by $\rho^{-n}f^{n*}p\eta$.  The seminorm is then
$$|\rho^{-n}f^{n*}p\eta|_T^2 =
\rho^{-2n}|p\eta|^2_{f_*^nT}=\rho^{-n}|p\eta|^2_{\rho^{-n}f_*^nT}.$$
Since $p\eta$ is smooth, and since $\rho^{-n}f^n_*T$ converges to 
$c\mu^-$, it follows
that $|p\eta|_{\rho^{-n}f_*^nT}$ converges to $|p\eta|_{c\mu^-}$. Thus
$\rho^{-n}f^{n*}p\eta$ converges to zero, and we conclude that
$\{g_n^+\eta\}$ is Cauchy.

Now observe that
$dd^cg^+_n\beta=\rho^{-n}f^{n*}\beta-\beta$, so 
$g^+_n\beta\ge-\beta$.  Further,
$\{g^+_n\beta\}$ is essentially monotone by the last assertion in 
Theorem \ref{prelim25}.
  Thus  $\{g^+_n\beta\}$ is essentially a
regularizing sequence, and it follows that
$g^+$ satisfies (\ref{tregular}) by Theorem \ref{nrg7}.

Now we have shown that $g^\pm$ satisfies (\ref{tregular}) with 
$T=\beta+\mu^\mp$.  If we
set $v=g^+$, $u=g^-$, and $T=\beta$, then it follows from Proposition 
\ref{prop32}  that
$g^+\in L^1(\mu^-\wedge\beta)$.
\qed
\begin{cor}\label{muwelldefined}
If (\ref{eq3}) holds, then the expression (\ref{defofmu}) defines a 
probability measure
$\mu$ such that $g^\pm\in L^1(\mu)$.
\end{cor}
\proof  This follows from Proposition \ref{prop32} with $v=u=g^+$ and 
$T=\mu^-$. \qed

\begin{cor}
\label{mu3}
For any
smooth closed $(1,1)$ forms $\eta_1,\eta_2$ on $X$, we have
$$\lim_{n,m\to\infty} \frac{f^{n*}\eta_1}{\rho^n}\wedge
                     \frac{f^m_*\eta_2}{\rho^m}  =
\mu\cdot\pair{\eta_1}{\theta^-}\pair{\eta_2}{\theta^+},
$$
in the weak* topology on Borel measures.
\end{cor}

\begin{thm}
\label{mu17}
$\mu$ does not charge points.
\end{thm}

\proof
By Corollaries \ref{prelim228} and \ref{lelongnumber}, there is no
point in $X$ at which both $\mu^+$ and $\mu^-$ have positive Lelong
number.  So given $x\in X$ we may assume that $\nu(\mu^-,x) = 0$.
Choose a local coordinate system such that $x=0$ and $U$ is the unit 
ball, and  choose a
smooth, compactly supported function $\chi:U\to [0,1]$ such that
$\chi(z) = 1$ for $\norm{z}$ small enough.  Set $\chi_r(z) = \chi(r^{-1}z)$ for
all $0<r<1$.
We have
\begin{eqnarray*}
\mu(x) & \leq & \int\chi_r\,\mu \leq
\left|\int \chi_r\,\omega^+\wedge\mu^-\right|
+
\left|\langle d\chi_r, d^c(g^+\wedge\mu^-) \rangle \right| \\
& \leq & C_1r^2 + C_2|\chi_r|_{\mu^-}
\leq C_1r^2 + \frac{C_2}{r^2}\int_{B_x(r)}\beta\wedge\mu^-.
\end{eqnarray*}
The last expression on the right converges to a multiple of 
$\nu(\mu^-,x) = 0$ as $r\to 0$.
\qed

Since $\mu$ does not charge $I(f)$, its pushforward
$f_*\mu$ is well-defined.

\begin{thm} \label{invariance}$\mu$ is $f$-invariant.
\end{thm}

\proof
Let $\psi:X\to\R$ be a continuous function.  Then by definition and
Corollary \ref{mu3}
$$
\int \psi\,f_*\mu = \int (\psi\circ f)\,\mu
   = \lim_{n\to\infty} \int (\psi\circ f)\,\frac{f^{n*}\beta}{\rho^n}\wedge
                       \frac{f^n_*\beta}{\rho^n}.
$$
The measure $f^{n*}\beta\wedge f^n_*\beta$ does not charge points because
there is no point in $X$ at which both factors are singular.  Both
factors are smooth at points outside a finite set, so this measure does not
charge curves either.  Thus
\begin{eqnarray*}
\int_X (\psi\circ f)\, \frac{f^{n*}\beta}{\rho^n}\wedge
                       \frac{f^n_*\beta}{\rho^n}
& = &
\int_{X-(\cC(f)\cup I(f))} \psi\circ f\, \frac{f^{n*}\beta}{\rho^n}\wedge
                       \frac{f^n_*\beta}{\rho^n} \\
& = &
\int_{X-(\cC(f^{-1})\cup I(f^{-1}))} \psi\,
                    \frac{f^{(n-1)*}\beta}{\rho^{n-1}}\wedge
                    \frac{f^{n+1}_*\beta}{\rho^{n+1}} \\
& = &
\int_X \psi\,
                    \frac{f^{(n-1)*}\beta}{\rho^{n-1}}\wedge
                    \frac{f^{n+1}_*\beta}{\rho^{n+1}} ,
\end{eqnarray*}
where second equality holds because
$f:X-\cC(f) - I(f)\to X-\cC(f^{-1})-I(f^{-1})$ is a biholomorphism.  If we let
$n\to\infty$, then the last expression on the right converges to
$\int\psi\mu$.  We conclude that
$f_*\mu =
\mu$.
\qed

\begin{cor}
\label{mu18}

$\mu$ does not charge $\cC(f)$ or $\cC(f^{-1})$.
\end{cor}

\proof
We already know that $\mu(I(f)) = 0$.  So by invariance,
$$
\mu(\cC(f)) = \mu(\cC(f)-I(f)) = \mu(I(f^{-1})) = 0.
$$
Similarly, $\mu(\cC(f^{-1})) = 0$.
\qed

\section{Mixing}
\label{mixing}
Let $T$ be a positive, closed (1,1) current, and let $J$ denote the
operator on 1-forms induced by the complex structure.  For a smooth
1-form
$\eta$, $\eta\wedge J\eta$ is a positive (1,1)-form, and we define the
$L^2(T)$ seminorm
$$||\eta||^2_{L^2(T)}:= \int \eta\wedge J\eta\wedge T.$$
If $\varphi\in L^1(T\wedge\beta)$, then we may define the quantity
$$|\varphi|^\natural_T:=\sup\{\int \varphi\wedge d\eta\wedge T: \eta{\rm\
  a\ smooth\ \hbox{1-form}\ with\ } ||\eta||_{L^2(T)}\le1\};$$
$|\cdot|_T^\natural$ is a seminorm on the space $\{\varphi\in
L^1(T\wedge\beta):|\varphi|_T^\natural<\infty\}$.
If $\varphi$ is smooth, then we may take $\eta=d^c\varphi$
and integrate by parts in the integral defining
$|\varphi|^\natural_T$ and apply Cauchy's inequality to see that
$|\varphi|^\natural_T =|\varphi|_T$.  Thus $|\cdot|_T^\natural$ extends
$|\cdot|_T$ from the space of smooth functions to a larger space.
\begin{lem}
\label{mixing1}  Suppose that (\ref{eq3}) holds.  If $\psi:X\to\R$ is a
smooth function, then
$$
|\psi\circ f^n|_{\mu^+}^\natural \leq \rho^{-n/2}|\psi|_{\mu^+}.
$$
\end{lem}

\proof
Recall that $f^n:X-(\cC(f^n)\cup I(f^n))\to X-(\cC(f^{-n})\cup I(f^{-n}))$ is
biholomorphic.  The function  $f^{n*}\psi$ is smooth and bounded on
$X-I(f^n)$, and $\beta\wedge\mu^+$ puts no mass on $I(f^n)$, so
$f^{n*}\psi\in L^1(\beta\wedge\mu^+)$.  Since $\beta\wedge\mu^+$ also puts no
mass on $\cC(f^n)$, it follows that if $\eta$ is a smooth 1-form, then
$$-\int(f^{n*}\psi)d\eta\wedge\mu^+ =-\int_{X-(\cC(f^n)\cup I(f^n))}
(f^{n*}\psi)d\eta\wedge\mu^+  = -\int_{X-(\cC(f^n)\cup I(f^n))}
f^{n*}(\psi{\mu^+\over\rho^n})\wedge d\eta.$$

Now let $\Gamma$ denote the (desingularized) graph of $f^n$, and let
$\pi_1,\pi_2:\Gamma\to X$ be the associated projections.  Thus
$f^n=\pi_2\pi_1^{-1}$, and  $f^{n*}=\pi_{1*}\pi_2^*$.  Now
$\pi_1:\Gamma-\cC(\pi_1)\to X-I(f^n)$ is a biholomorphism, so we may pull
this last integral back under $\pi_1^*$ to obtain an integral over $\Gamma$:
$$=-\int_{\Gamma-\cC(\pi_1)}\pi_2^*(\psi\wedge{\rho^{-n}\mu^+})\wedge
\pi_1^*(d\eta) =
-\int_{\Gamma-\cC(\pi_1)}\pi_2^*(\psi)\wedge\pi_1^*(d\eta)\wedge\pi_2^*
({\rho^{-n}\mu^+}).$$
Now $\pi_2^*\psi$ and $\pi_1^*(d\eta)$ are smooth on $\Gamma$, and
$\pi_2^*\psi$ and $\pi_1^*(d\eta)$ are smooth on $\Gamma$, and
$\pi_2^*({\rho^{-n}\mu^+})$ puts no mass on
$\cC(\pi_1)=\pi_1^{-1}(I(f^n))$, so we obtain
$$=-\int_{\Gamma} \pi_2^*(\psi)\wedge\pi_1^*(d\eta)\wedge
\pi_2^*\left({\rho^{-n}\mu^+}\right) =\int_{\Gamma}
\pi^*_2(d\psi)\wedge\pi_1^*(\eta)\wedge\pi_2^*\left({\rho^{-n}\mu^+}\right).$$
Applying the Schwarz inequality to this last term, we have
$$\left|\int(f^{n*}\psi)\wedge d\eta\wedge\mu^+\right|\le \left(\int_\Gamma
\pi^*_2(d\psi)\wedge\pi_2^*(d^c\psi)\wedge
\pi_2^*\left({\rho^{-n}\mu^+}\right)\right)^{1\over2}
\left(\int_\Gamma
\pi_1^*(\eta)\wedge J\pi^*_1(\eta)\wedge\pi_2^*\left({\rho^{-n}\mu^+}\right)
\right)^{1\over2}.$$
Again, since $\pi^*_2(d\psi)$ and $\pi_1^*(\eta)$ are smooth, there is no
change if we integrate over just $\Gamma-\cC(\pi_1)$ or $\Gamma-\cC(\pi_2)$.
Pushing the last integral forward under $\pi_{1*}$ and the first integral
forward under $\pi_{2*}$, we obtain
$$\left|\int(f^{n*}\psi)\wedge d\eta\wedge\mu^+\right|\le 
\left(\int_{X-I(f^{-n})}
d\psi\wedge d^c\psi\wedge{\rho^{-n}\mu^+}\right)^{1\over2}
\left(\int \eta\wedge J\eta\wedge
f^{n*}\left({\rho^{-n}\mu^+}\right)\right)^{1\over2}=$$
$$
=\rho^{-n/2}|\psi|_{\mu^+}||\eta||_{L^2(\mu^+)}.$$
Now we take the supremum over $\eta$ with $||\eta||_{L^2(\mu^+)}\le 1$, and the
Lemma follows.  \qed

\begin{thm}
\label{mixing4}
Let $\psi:X\to\R$ be a smooth function.  Then
$$
\lim_{n\to\infty} (\psi\circ f^n)\mu^+ = \mu^+\cdot\int\psi\,\mu.
$$
\end{thm}

\proof
We assume without loss of generality that $0\leq \psi \leq 1$ so that
$0 \leq (\psi\circ f^n)\mu^+ \leq \mu^+$.   Thus every subsequence of
$\{(\psi\circ f^n)\mu^+\}$ has a sub-subsequence converging to a positive
$(1,1)$ current $S$ dominated by $\mu^+$.  To see that $S$ is closed, let
$\eta$ be any smooth 1-form,
$$|\langle \eta,dS\rangle| = \lim_{n_j\to\infty}|\langle -d\eta,(\psi\circ
f^{n_j})\mu^+\rangle| \le ||\eta||_{L^2(\mu^+)}\,|\psi\circ 
f^n|^\natural_{\mu^+}.$$
The right hand side tends to zero by Lemma
\ref{mixing1}.  The current $\mu^+$ is  extremal in its cohomology class
[Gue, Theorem 5.1], so $S=c\mu^+$.  Thus it will suffice to show that the
cohomology class of $S$ is $c\theta^+$, with $c=\int\psi\mu$.

Let
$\pi_1,\pi_2:\Gamma\to X$ be as in the proof of Lemma \ref{mixing1}, so
$f^{n*}=\pi_{1*}\pi_2^*$.  Since $\beta$ is smooth and $\mu^+$ puts no mass on
$\cC(f^{-1})\cup\cC(f)$, we may pull integrals up from $X$ to 
$\Gamma$ and push them
back down as before:
\begin{eqnarray*}
\int_X (\psi\circ f^n)\, \mu^+\wedge\beta & = &
\int_X  \left( f^{n*}\psi\right)\wedge 
\left({\rho^{-n}}f^{n*}\mu^+\right) \wedge
\beta =
\int_\Gamma
(\pi_1^*\psi)\wedge \left({\rho^{-n}}{\pi_1^*\mu^+}\right)
  \wedge\pi_2^*\beta\\
& = &
\int_X 
\psi\wedge\pi_{1*}\left(\pi_2^*\beta\wedge\frac{\pi_1^*\mu^+}{\rho^n}\right) 
=
\int_X \psi\,\frac{f^n_*\beta}{\rho^n}\wedge\mu^+.
\end{eqnarray*}
Note that the second and fourth equalities depend on the fact that none of the
measures involved charge curves.  As
$n\to\infty$, the last integral converges to
$
\left(\int \psi\,\mu\right)\pair{\beta}{\theta^+}
$ by Corollary 4.5.
\qed

\begin{thm}
\label{mixing11}
The measure $\mu$ is mixing for $f$.
\end{thm}

\proof  We show that for any smooth functions $\varphi$ and $\psi$ we have
$$\lim_{n\to\infty}\int (\psi\circ f^n) \cdot \varphi\,\mu =
      \int \psi\,\mu\,\,\cdot\,\int\varphi\,\mu.$$
Let $\{g_j^-\}$ be a regularizing sequence for $g^-$.
By Theorem \ref{mixing4}, we have
\begin{equation}
\label{proof1}
\lim_{n\to\infty}
\int(\psi\circ f^n)\mu^+\wedge\varphi\wedge (\omega^-+dd^cg^-_j)=
\int \varphi\mu^+\wedge (\omega^-+dd^cg^-_j) \cdot\int\psi\,\mu.
\end{equation}
By Theorem \ref{mu2},
$\mu^+\wedge(\omega^-+dd^cg^-_j)$ converges to $\mu$ as $j\to\infty$.  Thus it
suffices to show that we may interchange the limits $n\to\infty$ and 
$j\to\infty$ in
the left hand integral in (\ref{proof1}).
For this, it suffices to show that
\begin{equation}
\label{equation2}
\lim_{j,k\to\infty}\int (\psi\circ f^n)\varphi \mu^+\wedge
dd^c(g^-_j-g^-_k)=0
\end{equation}
holds uniformly in $n$.
Integration by parts gives
$$\int(\psi\circ f^n)\varphi \mu^+\wedge
dd^c(g^-_j-g^-_k) =$$
$$= -\int \varphi \wedge d(\psi\circ f^n) \wedge
d^c(g^-_j-g^-_k)\wedge \mu^+ - \int(\psi\circ f^n)d\varphi \wedge
d^c(g^-_j-g^-_k)\wedge \mu^+ = I + II.
$$
The first integral is estimated by
$$|I| \le |\psi\circ f^n|^\natural_{\mu^+} \, ||\varphi
d^c(g^-_j-g^-_k)||_{L^2(\mu^+)}
\le  |\psi\circ f^n|^\natural_{\mu^+} \, ||\varphi||_{L^\infty} \,
|g_j^--g_k^-|_{\mu^+}.$$
This term converges to zero as $j,k\to\infty$ because $|\psi\circ
f^n|^\natural_{\mu^+} $ tends to zero by Lemma \ref{mixing1} and 
because $\{g^-_j\}$
is Cauchy for $|\cdot|_{\mu^+}$ by Theorem \ref{mu2}.
The second term is estimated by
$$|II| \le \left(\int (\psi\circ f^n)^2 d\varphi\wedge d^c\varphi\wedge
\mu^+\right)^{1/2} |g^-_j-g^-_k|_{\mu^+} \le ||\psi||_{L^\infty}\,
|\varphi|_{\mu^+}\, |g^-_j-g^-_k|_{\mu^+},
$$
which converges to zero because $\{g_j^-\}$ is Cauchy.  \qed

We can now strengthen Corollary \ref{mu18}.
\begin{cor}
\label{mixing5}
Every compact curve in $X$ has zero $\mu$-measure.
\end{cor}

\proof
Let $V$ be a compact irreducible curve.  By Corollary \ref{mu18}, we can
assume that
$V$ is not critical for any iterate of $f$.  If $V$ is not fixed by 
any iterate of
$f$, then $f^n(V)\cap f^m(V)$ is finite for any $n\neq m$.  Hence,
$$
\infty > \mu(X) \geq \mu(\bigcup_{n=0}^\infty f^n(V))
        = \sum_{n=0}^\infty \mu(f^n(V))
        = \sum_{n=0}^\infty \mu(V)
$$
which can happen only if $\mu(V) = 0$.

Now suppose that $V = f^n(V)$ is fixed by an iterate of $f$.
Since $\mu$ is mixing, we see that $V$ has either full or zero $\mu$-measure.
In the former case, it would follow that $f^n|_V$ lifts to an automorphism of
a compact Riemann surface (the desingularization of $V$) with a non-trivial
mixing invariant measure.  No such automorphism exists in dimension one,
so we must have $\mu(V) = 0$.
\qed

\begin{thm}
Any bimeromorphic map $h:X\to Y$ is defined at $\mu$ almost every 
point.  The pushforward
$h_*\mu$ is a probability measure, does not charge compact curves in $Y$, and
is invariant and mixing for $g\eqdef h^{-1}\circ f\circ h:Y\to Y$. 
If $g$ satisfies
conditions (\ref{eq3}), and $\mu^+_g$ and $\mu^-_g$ denote the associated
currents, then $h_*\mu = \mu^+_g\wedge\mu^-_g$.
\end{thm}

\proof  Let $\tilde h$ denote the restriction of $h$ to $X-(I(h)\cup \cC(h))$.
All but the last conclusion are immediate from Corollary 
\ref{mixing5} and the fact that
$\tilde h:X-(I(h)\cup \cC(h))\to Y-(I(h^{-1})\cup\cC(h^{-1}))$ is a 
biholomorphism.
Let us use the notation $\widetilde{\cdot}$ denote the restriction of 
a current or
measure on $Y$ to $Y-(I(h^{-1}) \cup \cC(h^{-1}))$.  With this 
notation, it follows
from  Corollary \ref{prelim6} that  $\tilde h^*\tilde\mu_g^+$ is equal to the
restriction of $\mu^+$ to
$X-(I(h)\cup
\cC(h))$.  Thus, since $h_*=(\tilde h^*)^{-1}$ on $X-(I(h)\cup 
\cC(h))$, we have
$
h_*\mu = \widetilde{h_*\mu}
= \widetilde{h_*\mu^+}\wedge\widetilde{h_*\mu^-} = \mu^+_g\wedge\mu^-_g.
$
\qed

\section{Lyapunov Exponents}
\label{lyapunov}
In this Section we will show
$\mu$ is hyperbolic (Theorem \ref{lyap4}).  In order to do this, we first
show that the Lyapunov exponents are finite $\mu$ a.e.  Next we give an
estimate on
$\mu$ of a ball with a certain mapping property (Proposition
\ref{lyap3}).  Then, using the machinery of Pesin Theory, we
convert this estimate into a proof of Theorem \ref{lyap4}.

\begin{prop}
\label{lyapunov1}
The quantities $\log^+\norm{Df}$, $\log^+\norm{(Df)^{-1}}$, and
$\log^+\norm{D^2 f}$ are $\mu$ integrable.
\end{prop}

\proof  It suffices to consider a neighborhood of a point $x\in I(f)$;
let $z$ be a local coordinate system such that $x$ corresponds to $z=0$.
By Proposition
\ref{derivative}, we have $\log^+||Df||\le C\left|
\log||z||\,\right|$.  Now if $v=\log||z||$, $u=g^+$, and $T=\mu^-$, we
may apply Theorems \ref{mu2} and \ref{nrg5} and Proposition \ref{prop32}
to conclude that
$v$, and thus
$\log^+||Df||$, is
$\mu$ integrable.  Likewise,  $\log^+\norm{D(f^{-1})}$ is integrable with
respect to
$\mu$.  Since $\mu$ is $f$-invariant, this gives the $\mu$-integrability
of
$\log^+\norm{(Df)^{-1}}$.

By the Cauchy estimates applied to the entries of $Df$ we have
$\log||D^2f||\le  C'\left|
\log||z||\,\right|$, so $\log||D^2f||$ is $\mu$-integrable.
\qed

Proposition \ref{lyapunov1} allows us to invoke Oseledec's Theorem
(see [KH, Theorem S.2.9]) to conclude:
\begin{prop} The limits
$$
\chi^+ = \lim_{n\to\infty} \frac{1}{n}\log||D f^n_x||,
\qquad
\chi^- = -\lim_{n\to\infty} \frac{1}{n}\log||D f^{-n}_x||
$$
exist and are finite at $\mu$ almost every point $x\in X$.    Since
$\mu$ is ergodic, these limits are a.e. independent of $x$.
\end{prop}

For
$v\in T_xX-\{0\}$, the Lyapunov exponent of
$f$ at
$x$ in the direction $v$ is defined as the limit
$$\chi(x,v)=\lim_{n\to\infty}{1\over n}\log||Df^n_xv||,$$
provided that this limit exists.
If $\chi^+\ne\chi^-$, it is a consequence of the
Oseledec Theorem is that there is an
$f$-invariant splitting
$T_xX=E^s_x\oplus E^u_x$ for $\mu$ a.e.\ $x$, and
$$\chi^\pm=\chi(x,v)=\lim_{n\to\infty}{1\over n}\log|Df_xv|,$$
for $v\in E^{u/s}_x-\{0\}$.

\begin{prop}
\label{lyap3}
There exist $C<\infty$ and $R>0$ such that if $f^n(\ball{x}{2r}) \subset
\ball{y}{R}$ for some $n\ge0$, then
$$
\mu(\ball{x}{r}) \leq C\rho^{-n/2}.
$$
\end{prop}

\proof
 From invariance of $\mu$ and Lemma \ref{mixing1}, we obtain
\begin{eqnarray*}
\left|\int \varphi\,\mu\right|
& = &
\left|\int (\varphi\circ f^{-n})\,\mu\right|
\leq
\left|\int (\varphi\circ f^{-n})\,\omega^+\wedge\mu^-\right|
+
\left|\int d(\varphi\circ f^{-n})\wedge d^cg^+\wedge\mu^-\right| \\
& \leq &
\frac{1}{\rho^n}\left|\int \varphi\wedge(f^{n*}\omega^+)\wedge\mu^-
\right| +
|g^+|_{\mu^-}|\varphi\circ f^{-n}|_{\mu^-} \\
& \leq &
\frac{1}{\rho^n}\left|\int \varphi\wedge(f^{n*}\omega^+)\wedge\mu^-
\right| +
\frac{C|\varphi|_{\mu^-}}{\rho^{n/2}}
\end{eqnarray*}
for any $\varphi\in C^\infty(X)$.  Now let us take $\varphi$ to be bounded
above and below by the characteristic functions for $\ball{x}{2r}$ and
$\ball{x}{r}$, respectively.  We may choose $\varphi$ such that
$\norm[C^1]{\varphi}^2,\norm[C^2]{\varphi}\leq Cr^{-2}$ for some constant
$C$ independent of $x$ and $r$.

Let us choose $R>0$ small enough that there is a local potential $u$ for
$\omega^+$ on $B_y(R)$.  Since $\omega^+$ is smooth, we may assume that
the $L^\infty$ norm of $u$ is bounded above independently of $y$.  We use
the local potential $dd^cu=\omega^+$ on the first right hand term in the
inequality above and integrate by parts twice to find:
\begin{eqnarray*}
\left|\int \varphi\,\mu\right| & \leq &
\frac{1}{\rho^n}\left|\int_{\ball{p}{2r}}
(u\circ f^n)\,dd^c\varphi\wedge\mu^- \right|
+
\frac{C}{\rho^{n/2}}
\left(\int_{\ball{p}{2r}} d\varphi\wedge d^c\varphi\wedge\mu^-\right)^{1/2} \\
& \leq &
\frac{C_1}{r^2\rho^n}\int_{\ball{p}{2r}} \beta\wedge\mu^-
+
\frac{C_2}{r\rho^{n/2}}\left(\int_{\ball{p}{2r}}
\beta\wedge\mu^-\right)^{1/2} \\
& \leq & \frac{C_1}{\rho^n} + \frac{C_2}{\rho^{n/2}}
\leq \frac{C}{\rho^{n/2}},
\end{eqnarray*}
where $C$ does not depend on $x$, $y$, or $r$.
\qed
\begin{thm}
\label{lyap4} The Lyapunov exponents satisfy:
$$
\chi^- \leq -{\log\rho\over 8} < 0 < {\log\rho\over8} \leq \chi^+.
$$
In particular, $\mu$ is hyperbolic of saddle type.
\end{thm}

\proof  Pesin Theory provides us with the following setup (see
[KH, Theorem S.3.1]).  For  any $\epsilon>0$, there exist
measurable radius and distortion functions $r, A:X\to\R^+$, and a
constant $c>0$ with the following properties.  For $\mu$ almost
every
$x\in X$ there is an embedding
$
\psi_x:\ball{\origin}{r(x)}\to X
$
such that:
\begin{itemize}
\item $\psi_x(\origin) = x$;
\item both $r$ and $A$ are `$\epsilon$-slowly varying,' i.e.,
$-\epsilon < \log A(fx)/A(x), \log r(fx)/r(x) < \epsilon;$
\item $c\, \dist(\psi_x(a),\psi_x(b)) \leq \dist(a,b)
      \leq A(x) \dist(\psi_x(a),\psi_x(b))$ for $a,b\in B_x(r(x))$;
\item if $f_x = \psi_{f(x)}^{-1}\circ f\circ \psi_x$, then
$
D_{\origin} f_x
$
is a diagonal matrix with diagonal entries
$e^{\chi_1},e^{\chi_2}$ satisfying
$$
|\mathrm{Re}\,\chi_1 - \chi^+|, |\mathrm{Re}\,\chi_2 - \chi^-| < \epsilon;
$$
\item $\norm[C^1]{f_x - D_{\origin}f_x} \leq \epsilon$ on the domain
       of $f_x$.
\end{itemize}
With this notation, we set
$r_n(x) = e^{-n(\chi^+ + 3\epsilon)}cr(x)/A(x)$.  By the properties
above, it follows that for
$1\leq j \leq n$ we have
$$
f^j|_{\ball{x}{r_n(x)}}
= \psi_{f^j(x)}\circ f_{f^{j-1}(x)}\circ\dots\circ f_x\circ \psi_x^{-1}.
$$
That is, by keeping track of the diameters of the successive images of
$\ball{x}{r_n(x)}$ one sees that each stage of the composition on the
right side is well-defined and that, moreover,
$$
f^n\ball{x}{r_n(x)} \subset \ball{f^n(x)}{r(f^n(x))}
$$
Therefore Lemma \ref{lyap3} gives us the bound
$$
\mu(\ball{x}{r_n(x)}) \leq C\rho^{-n/2}.
$$
Now Lusin's Theorem provides us with a compact
subset $K\subset\supp\mu$ such that $\mu(K) > 1/2$ and on which
$r$ and $A$ vary continuously.  Thus for all $x\in K$ the radius of
$B_n(x)$ is bounded below by $C e^{-n(\chi^+ + 3\epsilon)}$.  We can
therefore choose $m \leq C e^{4n(\chi^+ + 3\epsilon)}$ points
$x_1,\dots, x_m\in K$ such that
$K\subset\bigcup_{j=1}^m \ball{x_j}{r_n(x_j)}$.
Using this cover, we estimate
$$
1/2 < \mu(K) \leq \sum_{j=1}^m \mu(\ball{x_j}{r_n(x_j)})
              \leq  C e^{4n(\chi^+ + 3\epsilon)}\rho^{-n/2}.
$$
Letting $n$ tend to $\infty$ and then $\epsilon$ to $0$ yields
$$
1 < e^{4\chi^+}\rho^{-1/2},
$$
which yields the estimate of Theorem \ref{lyap4}.
\qed

Using
Theorem \ref{lyap4} and the fact that $\mu$ is mixing, one can
apply a contraction mapping argument, (see, for example [Dil2, \S8]), to
obtain:

\begin{cor} $\supp\mu$ is contained in the closure of the saddle
periodic points of $\mu$.
\end{cor}

\bigskip

\centerline{\bf References}
\medskip

[BPV]  W. Barth, C. Peters, and A. van de Ven, {\sl Compact Complex
Surfaces}, Springer-Verlag, Berlin, 1984.
\smallskip

[BS] Eric Bedford and John Smillie, Polynomial diffeomorphisms of $ C\sp
2$: currents, equilibrium measure and hyperbolicity. Invent. Math. 103
(1991), no. 1, 69--99.
\smallskip

[BT]  E. Bedford and B.A. Taylor, Variational properties of the complex
Monge-Amp\`ere equation I. Dirichlet principle, Duke Math. J. 45 (1978),
375--403.
\smallskip

[B]  Z. Blocki, On the definition of the Monge-Amp\`ere operator in ${\bf
C}^2$, preprint.
\smallskip

[BD]  Jean-Yves Briend and Julien Duval, Exposants de Liapounoff et
distribution des points p\'eriodiques d'un endomorphisme de ${\bf C}P^k$,
Acta Math. 182 (1999), 143--157.
\smallskip

[C] Serge Cantat,  Dynamique des automorphismes des surfaces $K3$.
(French) [Dynamics of $K3$-surface automorphisms] Acta Math. 187 (2001),
no. 1, 1--57.
\smallskip

[Dem1]  Jean-Pierre Demailly, Regularization of closed positive currents
and intersection theory.  J. Algebraic Geom., 1 (1992),
361--409. \smallskip

[Dem2]  Jean-Pierre Demailly, Monge-Amp\`ere operators, Lelong numbers
and intersection theory.  In {\sl Complex Analysis and Geometry}, pages
115--193.  Plenum, New York, 1993.
\smallskip

[Dil1]  Jeffrey Diller, Dynamics of birational maps of ${\bf P}^2$,
Indiana U. Math. J. 45 (1996), 721--772.\smallskip

[Dil2] Jeffrey Diller, Invariant Measure and Lyapunov exponents for
separating birational maps, Comment. Math. Helv., 76 (2001), no. 4,
754--780.
\smallskip

[Dil3] Jeffrey Diller, Birational maps, positive currents, and dynamics.
Michigan Math. J. 46 (1999), no. 2, 361--375. \smallskip

[DF]  Jeffrey Diller and Charles Favre,  Dynamics of bimeromorphic maps of
surfaces. Amer. J. Math. 123 (2001), no. 6, 1135--1169.
\smallskip

[Fav1] Charles Favre, Note on pull-back and Lelong number of currents,
Bull. Soc. Math. France 127 (1999), 445--458.
\smallskip

[Fav2] Charles Favre, Multiplicity of holomorphic functions, Math. Ann.
316 (2000), 355--378.
\smallskip

[Fav3] Charles Favre,  Points p\'eriodiques d'applications birationnelles
de ${\bf P}^2$. (French) [Periodic points of birational mappings of
  ${\bf P}^2$] Ann. Inst. Fourier (Grenoble) 48 (1998), no. 4,
999--1023.
\smallskip

[FS1] John-Erik Forn\ae ss and Nessim Sibony, Complex H\'enon mappings in
${\bf C}^2$ and Fatou-Bieberbach domains. Duke Math. J. 65 (1992), no. 2,
345--380.
\smallskip

[FS2] John-Erik Forn\ae ss and Nessim Sibony, Complex dynamics in higher
dimension. I. Complex analytic methods in dynamical systems (Rio de
Janeiro, 1992). Ast\'erisque No. 222 (1994), 5, 201--231.
\smallskip

[GH]  Phillip Griffiths and Joseph Harris, {\sl Principles of Algebraic
Geometry}, John Wiley \&\ Sons Inc., New York, 1994.
\smallskip

[Gue]  Vincent Guedj, Dynamics of polynomial mappings of ${\bf C}^2$,
Amer. J. Math. 124 (2002), no. 1, 75--106.
\smallskip

[KH]  Anatole Katok and Boris Hasselblatt, {\sl Introduction to the
Modern Theory of Dynamical Systems}. Cambridge University Press, 1995.
\smallskip

[Lam] A. Lamari, Le c\^one K\"ahl\'erien d'une surface, J. Math. Pures
Appl. (9) 78 (1999), 249--263.
\smallskip

\end{document}